\documentclass[a4paper,reqno]{amsart}

\usepackage[british]{babel}
\usepackage[T1]{fontenc}
\usepackage{amssymb,amsthm,bm}
\usepackage{mathtools}
\usepackage{mathrsfs}
\usepackage{nicefrac}
\mathtoolsset{mathic}
\allowdisplaybreaks
\usepackage{courier,mathptmx,stmaryrd}
\usepackage{longtable}
\usepackage[numbers,sort&compress]{natbib}
\usepackage{hyperref,url}
\hypersetup{  
bookmarksnumbered
}
\hypersetup{  
breaklinks=false,
pdfborderstyle={/S/U/W 0},  
citebordercolor=.235 .702 .443,
urlbordercolor=.255 .412 .882,
linkbordercolor=.804 .149 .19,
}
\hypersetup{ 
pdfauthor={Gert de Cooman, Arthur Van Camp and Jasper De Bock},
pdftitle={The logic behind desirable sets of things, and its filter representation},
pdfkeywords={desirability, desirable sets of things, conservative inference, propositional logic, filter, prime filter, principal filter}
}
\usepackage[all]{hypcap} 
\usepackage{xcolor}
\usepackage{enumitem}
\usepackage{tikz}
\usetikzlibrary{decorations.text}
\usetikzlibrary{shapes}

\newif\ifbibincluded
\bibincludedfalse

\newtheorem{theorem}{Theorem}
\newtheorem{proposition}[theorem]{Proposition}
\newtheorem{lemma}[theorem]{Lemma}
\newtheorem{corollary}[theorem]{Corollary}
\theoremstyle{definition}
\newtheorem{definition}{Definition}
\theoremstyle{remark}

\newtheorem{example}{Example}

\DeclarePairedDelimiter{\abs}{\vert}{\vert}

\DeclarePairedDelimiter{\ceil}{\lceil}{\rceil}
\DeclarePairedDelimiter{\group}{(}{)}
\DeclarePairedDelimiter{\sqgroup}{[}{]}

\newcommand{\set}[2][]{#1\{#2#1\}}
\newcommand{\cset}[3][]{\set[#1]{#2\colon#3}}

\newcommand{\naturals}{\mathbb{N}}
\newcommand{\naturalswithzero}{{\mathbb{N}_0}}
\newcommand{\integers}{\mathbb{Z}}
\newcommand{\reals}{\mathbb{R}}

\newcommand{\rationals}{\mathbb{Q}}

\newcommand{\posspace}{\mathscr{X}}
\newcommand{\outcomes}{\{0,1\}}

\newcommand{\gamble}{f}
\newcommand{\gambles}{\reals^\posspace}
\newcommand{\frcsts}{\mathscr{I}}
\newcommand{\ratfrcsts}{\smash{\mathscr{I}_{r}}}
\newcommand{\frcst}{I}
\newcommand{\ratfrcst}{I_r}
\newcommand{\frcstseqs}{{\mathscr{I}_{r}^{\naturals}}}
\newcommand{\frcstseq}{\iota}
\newcommand{\frcstsits}{{\mathscr{I}_{r}^{*}}}
\newcommand{\frcstsit}{i}
\newcommand{\frcstsystems}{\Phi}
\newcommand{\ratfrcstsystems}{\Phi_r}
\newcommand{\frcstsystem}{\varphi}
\newcommand{\ratfrcstsystem}{\varphi_r}
\newcommand{\lfrcstsystem}{\underline{\frcstsystem}}
\newcommand{\ufrcstsystem}{\overline{\frcstsystem}}
\newcommand{\lratfrcstsystem}{\underline{\frcstsystem}_r}
\newcommand{\uratfrcstsystem}{\overline{\frcstsystem}_r}
\newcommand{\sits}{\posspace^{*}}
\newcommand{\sit}{w}
\newcommand{\altsit}{\overline{\sit}}
\newcommand{\pths}{\posspace^{\naturals}}
\newcommand{\pth}{\omega}

\newcommand{\preqpths}{\group{\smash{\ratfrcsts\times\posspace}}^{\naturals}}
\newcommand{\realpreqpths}{\group{\frcsts\times\posspace}^{\naturals}}
\newcommand{\preqpth}{\upsilon}
\newcommand{\preqsits}{\group{\smash{\ratfrcsts\times\posspace}}^{*}}
\newcommand{\preqsit}{v}

\newcommand{\preqpthrep}{\group{\frcstseq,\pth}}
\newcommand{\preqsitrep}{\group{\frcstsit,\sit}}

\newcommand{\supermartins}[1][\frcstsystem]{\overline{\mathbb{M}}(#1)}
\newcommand{\supermartin}{M}

\newcommand{\tests}[1][\frcstsystem]{\overline{\mathbb{T}}(#1)}
\newcommand{\test}{T}
\newcommand{\superfarths}{\overline{\mathbb{F}}}
\newcommand{\superfarth}{F}
\newcommand{\unifarth}{U}
\newcommand{\selection}{S}

\newcommand{\ex}{E}
\newcommand{\lex}{\underline{\ex}}
\newcommand{\uex}{\overline{\ex}}

\newcommand{\uexfrcstsystem}{\overline{\ex}_{\frcstsystem(\sit)}}

\newcommand{\uexratfrcstsystem}{\overline{\ex}_{\ratfrcstsystem(\sit)}}

\newcommand{\init}{\square}
\newcommand{\pthat}[1]{\pth_{#1}}
\newcommand{\pthatn}{\pth_{n}}
\newcommand{\pthatnplus}{\pth_{n+1}}

\newcommand{\pthatkplus}{\pth_{k+1}}
\newcommand{\pthto}[1]{\pth_{1:#1}}
\newcommand{\pthton}{\pth_{1:n}}

\newcommand{\andoutcome}{\,\cdot}
\newcommand{\sitto}[1]{\sit_{1:#1}}
\newcommand{\sitat}[1]{\sit_{#1}}
\newcommand{\altsitto}[1]{\altsit_{1:#1}}
\newcommand{\altsitat}[1]{\altsit_{#1}}

\newcommand{\preqpthto}[1]{\preqpth_{1:#1}}
\newcommand{\preqpthton}{\preqpth_{1:n}}
\newcommand{\preqpthtok}{\preqpth_{1:k}}

\newcommand{\preqsitto}[1]{\preqsit_{1:#1}}

\newcommand{\frcstsitat}[1]{\frcstsit_{#1}}
\newcommand{\frcstseqat}[1]{\frcstseq_{#1}}
\newcommand{\frcstseqatnplus}{\frcstseqat{n+1}}
\newcommand{\frcstseqatkplus}{\frcstseqat{k+1}}
\newcommand{\frcstseqto}[1]{\frcstseq_{1:#1}}
\newcommand{\frcstseqton}{\frcstseqto{n}}
\newcommand{\frcstseqtonplus}{\frcstseqto{n+1}}

\newcommand{\ifandonlyif}{\Leftrightarrow}

\newcommand{\lsc}{\nearrow}
\newcommand{\bolleke}{\vcenter{\hbox{\scalebox{0.7}{\(\bullet\)}}}}

\usepackage{nicefrac}
\usepackage{aligned-overset}
\usepackage{adjustbox}
\usepackage{listings}
\newcommand{\sectionref}[1]{Section~\ref{#1}}
\newcommand{\figureref}[1]{Figure~\ref{#1}}
\newcommand{\definitionref}[1]{Definition~\ref{#1}}
\newcommand{\lemmaref}[1]{Lemma~\ref{#1}}
\newcommand{\equationref}[1]{Equation~\eqref{#1}}
\newcommand{\corollaryref}[1]{Corollary~\ref{#1}}
\newcommand{\theoremref}[1]{Theorem~\ref{#1}}
\newcommand{\propositionref}[1]{Proposition~\ref{#1}}
\newcommand{\exampleref}[1]{Example~\ref{#1}}

\title{
  Imprecision in Martingale-Theoretic Prequential Randomness 
}

\author{Floris Persiau \and Gert de Cooman}
\address{Ghent University, Foundations Lab, Technologiepark--Zwijnaarde 125, Ghent, Belgium}
\title{Imprecision in Martingale-Theoretic Prequential Randomness}

\begin{document}
\begin{abstract}
In a prequential approach to algorithmic randomness, probabilities for the next outcome can be forecast `on the fly' without the need for fully specifying a probability measure on all possible sequences of outcomes, as is the case in the more standard approach.
We take the first steps in allowing for probability intervals instead of precise probabilities in this prequential approach, based on ideas from our earlier imprecise-probabilistic and martingale-theoretic account of algorithmic randomness.
We define what it means for an infinite sequence~\((I_1,x_1,I_2,x_2,\dots)\) of successive interval forecasts~\(I_k\) and subsequent binary outcomes~\(x_k\) to be random.
We compare the resulting prequential randomness notion with the more standard one, and investigate where both randomness notions coincide, as well as where their properties correspond.
\end{abstract}
\maketitle

\section{Introduction}\label{sec:intro}
Consider an infinite sequence~\(X_1,X_2,X_3,\dots\) of binary variables~\(X_k\in\outcomes\).
In classical probability theory, uncertainty about the possible outcomes~\(X_k\), with~\(k\in\naturals\),\footnote{\(\naturals\) denotes the set of natural numbers, and \(\naturalswithzero\coloneqq\naturals\cup\set{0}\) the set of non-negative integers.\footnotemark}\footnotetext{A real~\(x\in\reals\) is called \emph{positive} if \(x>0\) and \emph{non-negative} if \(x\geq0\).} is typically represented by a probability measure on the elements of the Borel algebra over the set~\(\outcomes^\naturals\) generated by binary strings.
This is equivalent to specifying, for every possible binary string~\(\sit\in\outcomes^n\), for any~\(n\in\naturalswithzero\), the probability~\(p_\sit\in\sqgroup{0,1}\) that the variable~\(X_{n+1}\) equals~\(1\), given the observation of the binary string~\(\sit\); we will call this a \emph{standard} approach to (im)precise probability theory.
This system of (conditional) probabilities is also called a \emph{forecasting system}.
In a number of papers \cite{Dawid1997,Dawid1999}, Dawid and Vovk question whether it is always natural or even possible for a subject to specify such a probability measure.
Consider for example a weather forecaster who provides a daily probability for rain in the next \(24\) hours.
His forecasts are based on the rain history he has actually observed (as well as other information), and he isn't required to provide forecasts for all rain histories that might have been and might be.
Dawid and Vovk provide a practical way out of this conundrum by putting forward the so-called \emph{prequential forecasting} framework.

In Refs.~\cite{Shernov2008,Vovk2010}, the prequential forecasting idea is applied to algorithmic randomness.
Instead of defining what it means for an infinite binary sequence of outcomes, such as \((0,1,1,0,0,1,0,\dots)\), to be random for a forecasting system, they come up with randomness notions that consider the randomness of an infinite binary sequence of outcomes only with respect to the probabilities that are forecast along the sequence, that is, they define what it means for an infinite sequence~\((p_1,x_1,p_2,x_2,\dots,p_k,x_k,\dots)\) of probability forecasts~\(p_k\in\sqgroup{0,1}\) and subsequent outcomes~\(x_k\in\outcomes\) to be random.
They do this using a measure-theoretic as well as a martingale-theoretic approach; an infinite sequence is regarded as measure-theoretically random if there is no computable way to specify a set of measure zero containing this sequence, whereas an infinite sequence is regarded as martingale-theoretically random if there is no computable way to get arbitrary rich by betting on its elements \cite{Vovk2010}.

Here, we build upon the martingale-theoretic prequential approach to randomness by extending the probability forecasts~\(p_k\in\sqgroup{0,1}\) to so-called interval forecasts~\(\frcst_k\subseteq\sqgroup{0,1}\), and in doing so extend the range of applicability of both their and our own work \cite{CoomanBock2021}: we define what it means to be random for an infinite sequence~\((\frcst_1,x_1,\frcst_2,x_2,\dots)\) of interval forecasts~\(\frcst_k\) and subsequent outcomes~\(x_k\), and compare the resulting prequential imprecise-probabilistic randomness definition to our previously introduced \emph{standard} imprecise-probabilistic generalisation of Martin-Löf randomness \cite{CoomanBock2021}.

We structure our exposition as follows.
In \sectionref{sec:local:models}, we formally introduce interval forecasts and equip them with an interpretation in terms of betting games.
This allows us to discuss, in \sectionref{sec:game}, the basic ideas behind our \emph{standard} imprecise-probabilistic martingale-theoretic approach to randomness \cite{CoomanBock2021}.
After summarising some computability principles and properties in \sectionref{sec:comp}, we have enough mathematical equipment to formally introduce our \emph{prequential} imprecise-probabilistic notion of randomness in \sectionref{sec:randomness}; our terminology will follow Refs.~\cite{Dawid1999,Vovk2010}.
We compare the definitions and properties of our standard and prequential randomness notions in \sectionref{sec:comparison}.

We have gathered more technical results and proofs in Appendix~\ref{sec:appendix}.

\section{Imprecise Uncertainty Models}\label{sec:local:models}
Consider the binary \emph{sample space}~\(\posspace\coloneqq\outcomes\) and a \emph{variable}~\(X\) that may assume values in~\(\posspace\).
To describe a subject's uncertainty about the unknown value of~\(X\), we'll not only allow for precise probabilities~\(p\in\sqgroup{0,1}\), but also for more general closed \emph{interval forecasts}~\(\frcst\subseteq\sqgroup{0,1}\); the set of all closed interval forecasts is denoted by~\(\frcsts\).
Every interval forecast~\(\frcst\in\frcsts\) could be interpreted as a set of probabilities~\(p\in\frcst\) that a subject finds plausible for describing his probability that \(X=1\).
To obtain a different interpretation in terms of bets---which is the one we use here---we associate with every interval forecast~\(\frcst\in\frcsts\) the \emph{upper} and \emph{lower expectation} operators~\(\uex_\frcst,\lex_\frcst\colon\reals^\posspace\to\reals\) that associate with every so-called \emph{gamble}~\(\gamble\in\reals^\posspace\) an upper expectation
\begin{align}
\uex_\frcst(\gamble)&\coloneqq\max_{p\in\frcst}\sqgroup[\big]{p\gamble(1)+(1-p)\gamble(0)}\label{def:uex}
\shortintertext{and a lower expectation}
\lex_\frcst(\gamble)&\coloneqq\min_{p\in\frcst}\sqgroup[\big]{p\gamble(1)+(1-p)\gamble(0)}
=-\uex_\frcst(-\gamble).\label{def:lex}
\end{align}
These numbers are interpreted as a subject's lowest acceptable selling price and largest acceptable buying price, respectively, for the uncertain pay-off~\(\gamble(X)\).
This implies that our subject is willing to accept the uncertain pay-off~\(\gamble(X)-p\) for any buying price~\(p\leq\lex_\frcst(\gamble)\), and is willing to accept the uncertain pay-off~\(q-f(X)\) for any selling price~\(q\geq\uex_\frcst(\gamble)\); the collection of our subject's accepted gambles corresponds to those gambles~\(\gamble\in\gambles\) for which \(\lex_\frcst(\gamble)\geq 0\), that is, all gambles for which he expects a non-negative gain with respect to every probability~\(p\in\frcst\).
Vice versa, from the perspective of an opponent, our subject is willing to give away those gambles~\(\gamble\) for which \(\uex_\frcst(\gamble)\leq 0\), that is, all gambles for which he expects a non-negative loss for every probability~\(p\in\frcst\).
Our subject is indeterminate about accepting or giving away a gamble~\(\gamble\) when \(\lex_\frcst(\gamble)<0<\uex_\frcst(\gamble)\); this is illustrated in \figureref{fig:desirable:gambles}.
In what follows, we'll make extensive use of the upper expectation operator~\(\uex_\frcst\) and a number of its properties \cite{walley1991}:

\begin{proposition}
Consider any interval forecast~\(\frcst\in\frcsts\).
Then for all gambles~\(\gamble,g\in\reals^\posspace\), all sequences of gambles~\(\gamble_n\), \(n\in\naturalswithzero\) and all~\(\mu,\lambda\in\reals\) with~\(\lambda\geq 0\):
\begin{enumerate}[label=\upshape{C}\arabic*.,ref=\upshape{C}\arabic*,leftmargin=*,itemsep=0pt]
\item\label{axiom:coherence:bounds} \(\min\gamble\leq\uex_\frcst(\gamble)\leq\max\gamble\);\hfill{\upshape[boundedness]}
\item\label{axiom:coherence:homogeneity} \(\uex_\frcst(\lambda\gamble)=\lambda\uex_\frcst(\gamble)\);\hfill{\upshape[non-negative homogeneity]}
\item\label{axiom:coherence:subsupadditivity} \(\uex_\frcst(\gamble+g)\leq\uex_\frcst(\gamble)+\uex_\frcst(g)\);\hfill{\upshape[subadditivity]}
\item\label{axiom:coherence:constantadditivity} \(\uex_I(\gamble+\mu)=\uex_\frcst(\gamble)+\mu\);\hfill{\upshape[constant additivity]}
\item\label{axiom:coherence:monotonicity} if \(f\leq g\) then \(\uex_\frcst(\gamble)\leq\uex_\frcst(g)\);\hfill{\upshape[monotonicity]}
\item\label{axiom:coherence:uniform:convergence} if \(\lim_{n\to\infty}\gamble_n=\gamble\) then \(\lim_{n\to\infty}\uex_\frcst(\gamble_n)=\uex_\frcst(f)\). \hfill{\upshape[pointwise convergence]}
\end{enumerate}
\end{proposition}

\begin{figure}[h]
\centering
\begin{tikzpicture}[scale=1.4]

\fill[opacity=0.5,red!80!black] (0,0) -- (-5/2*9/10,5/6*9/10) -- (-5/2*9/10,-5/2*9/10) -- (5/6*9/10,-5/2*9/10) --cycle;
\fill[opacity=0.5,green!60!black] (0,0) -- (+5/2*9/10,-5/6*9/10) -- (5/2*9/10,5/2*9/10) -- (-5/6*9/10,5/2*9/10) --cycle;

\draw[->] (-2.25,0) -- (2.25,0) node[below right] {\(f(1)\)};
\draw[->] (0,-2.25) -- (0,2.25) node[above left] {\(f(0)\)};

\draw[very thick,red!80!black] (0,0) -- (-5/2*9/10,5/6*9/10);
\draw[very thick,red!80!black] (0,0) -- (5/6*9/10,-5/2*9/10);
\draw[very thick,green!60!black] (0,0) -- (+5/2*9/10,-5/6*9/10);
\draw[very thick,green!60!black] (0,0) -- (-5/6*9/10,+5/2*9/10);

\foreach \x in {-2,-1,1,2}
    \draw (\x,2pt) -- (\x,-2pt)
	node[anchor=north] {\small\(\x\)};
\foreach \y in {-2,-1,1,2}
    \draw (2pt,\y) -- (-2pt,\y) 
    node[anchor=east] {\small\(\y\)};
\node[anchor=north east] at (-2pt, -2pt)  {\small\(0\)};

\draw[red!80!black,line width=0.2pt,fill=red!80!black,x=1pt,y=1pt] ({-1/sqrt(2)},{1/sqrt(2)}) arc (135:315:1pt);
\draw[green!60!black,line width=0.2pt,fill=green!60!black,x=1pt,y=1pt] ({1/sqrt(2)},{-1/sqrt(2)}) arc (-45:135:1pt);

\node[anchor=west,fill=white,rounded corners,opacity=60] at (0.35, 1)  {\large\(\lex_\frcst(\gamble)\geq0\)};
\node[anchor=east,fill=white,rounded corners,opacity=60] at (-0.55, -1)  {\large\(\uex_\frcst(\gamble)\leq0\)};
\end{tikzpicture}
\caption{Let \(\frcst=\sqgroup{\nicefrac{1}{4},\nicefrac{3}{4}}\). The red, green and white regions depict the gambles~\(\gamble\in\gambles\) for which \(\uex_\frcst(\gamble)\leq0\), \(\lex_\frcst(\gamble)\geq0\) and \(\lex_\frcst(\gamble)<0<\uex_\frcst(\gamble)\), respectively.}
\label{fig:desirable:gambles}
\end{figure}

\section{Sequential and Prequential Games}\label{sec:game}
To put interval forecasts into practice, consider Frank Deboosere---a famous Belgian weatherman---whose daily job consists in making good forecasts about whether the sun will or won't shine on the next day.
This corresponds to a binary option space; we write \(1\) for a sunny day and \(0\) for a non-sunny day.
We formalise his forecasting task in the following forecasting protocol:\vspace{8pt}

\begin{adjustbox}{width=\columnwidth/3*2,keepaspectratio}
\begin{lstlisting}[mathescape]
FOR $n=1,2,3,\ldots:$
$\quad$ Forecaster Frank announces $\frcst_n\in\frcsts$.
$\quad$ Reality announces $x_{n}\in\posspace$.$\quad\quad\quad\quad\quad\quad\quad$
$\newline$
\end{lstlisting}
\end{adjustbox}

\noindent Intuitively, at each step~\(n\in\naturals\) in the protocol, \(\frcst_n\) expresses Frank's beliefs about~\(X_n=1\) after observing the \emph{outcomes}~\(\group{x_1,\dots,x_{n-1}}\).
Clearly, Frank can do a good or a bad forecasting job.
For example, if he forecasts~\(1\) at every time step, but it rains every day, then we might be inclined to say he's doing a bad job.
But if he forecasts~\(\nicefrac{1}{2}\) at every time step and it rains half of the time, then we could say he's doing a good job.
This brings us to the central question in this paper: when will we say that Frank makes good predictions, or more technically speaking, that his forecasts~\((\frcst_1,\dots,\frcst_n,\dots)\) are well-calibrated with the outcomes~\((x_1,\dots,x_n,\dots)\)?
The field of algorithmic randomness tries to answer this question by defining what it means for an infinite sequence~\((\frcst_1,x_1,\dots,\frcst_n,x_n,\dots)\) of forecasts~\(\frcst_n\) and subsequent outcomes~\(x_n\) to be `random'.

\subsection{The Standard Approach}
Before giving a first (standard) answer to this randomness question, we need some notation.

An infinite sequence of outcomes~\(\group{x_1,x_2,\dots,x_n,\dots}\in\pths\) is called a \emph{path} and is generically denoted by~\(\pth\).
A finite sequence of outcomes~\(\group{x_1, x_2,\dots,x_n}\in\sits\coloneqq\bigcup_{k\in\naturalswithzero}\posspace^k\) is called a \emph{situation}, is generically denoted by~\(\sit\), and has length~\(\abs{\sit}=n\).
For any~\(k\in\naturalswithzero\), we use the notations \(\pthto{k}=\group{x_1,x_2,\dots,x_k}\) and \(\pthat{k}=x_k\), and similarly for situations~\(\sit\in\sits\) with~\(k\leq\abs{\sit}\).
The empty situation~\(\sitto{0}=\group{}\) is also denoted by~\(\init\).

On the standard approach, it's assumed that Forecaster Frank's forecasts in the above protocol can be derived from a so-called \emph{forecasting system}: Frank not only has to specify forecasts~\(I_n\coloneqq I_{(x_1,\dots,x_{n-1})}\) to express his beliefs about~\(X_{n}=1\) \emph{after} observing the actual outcomes~\((x_1,\dots,x_{n-1})\), but he also has to specify forecasts~\(\frcst_\sit\in\frcsts\) for all possible situations~\(\sit\in\sits\) that can in principle occur.

\begin{definition}
A \emph{forecasting system} is a map~\(\frcstsystem\colon\sits\to\frcsts\) that associates an interval forecast~\(\frcstsystem(\sit)\in\frcsts\) with every situation~\(\sit\) in the event tree~\(\sits\).
With any forecasting system~\(\frcstsystem\) we can associate two real maps~\(\lfrcstsystem\) and \(\ufrcstsystem\), defined by~\(\lfrcstsystem(\sit)\coloneqq\min\frcstsystem(\sit)\) and \(\ufrcstsystem(\sit)\coloneqq\max\frcstsystem(\sit)\) for all~\(\sit\in\sits\).
We denote the set of all forecasting systems by~\(\frcstsystems\).
A forecasting system~\(\frcstsystem\in\frcstsystems\) is called \emph{non-degenerate} if \(\lfrcstsystem(\sit)<1\) and \(0<\ufrcstsystem(\sit)\) for all~\(\sit\in\sits\).
A forecasting system~\(\frcstsystem\in\frcstsystems\) is called \emph{more conservative than} a forecasting system~\(\frcstsystem'\in\frcstsystems\) if \(\frcstsystem'(\sit)\subseteq\frcstsystem(\sit)\) for all~\(\sit\in\sits\).
\end{definition}
\noindent We see that, in this context, it's more natural to talk about the randomness of a path~\(\smash{\pth\in\pths}\) for a forecasting system~\(\frcstsystem\), rather than for a sequence of forecasts~\((\frcst_1,\dots,\frcst_n,\dots)\).

To answer the randomness question, Frank's colleague Sabine Hagedoren---who is a famous Belgian weatherwoman and whom we'll also call Sceptic, because that will be her role---tests the correspondence between Frank's forecasting system~\(\frcstsystem\) and Reality's outcomes.
She does so by engaging in a betting game.
We'll assume that she starts with unit capital.
In every situation~\(\sit\in\sits\), she then selects a gamble~\(\gamble_\sit\in\gambles\) that's made available to her by Forecaster Frank's specification of the interval forecast~\(\frcstsystem(\sit)\), that is, she selects an uncertain change of capital~\(\gamble_\sit\in\gambles\) for which \(\uex_{\frcstsystem(\sit)}(\gamble_\sit)\leq0\).
Furthermore, we'll prohibit Sabine from borrowing money, which means that her capital can't become negative.
If Frank does a good forecasting job, Sabine shouldn't be able to tremendously increase her capital in the long run.
We'll then call a path~\(\pth\in\pths\) \emph{random} for a forecasting system~\(\frcstsystem\in\frcstsystems\) if Sabine can't come up with a(n effectively implementable) betting strategy that makes her rich without bounds along~\(\pth\); her betting strategies are formalised in the following definition, where `\(\andoutcome\)' functions as a placeholder for the possible outcomes~\(x\in\posspace\).

\begin{definition}
A real-valued map~\(\supermartin\colon\sits\to\reals\) is called a \emph{supermartingale} for~\(\frcstsystem\) if \allowbreak\(\uex_{\frcstsystem(\sit)}(\supermartin(\sit\andoutcome))\leq \supermartin(\sit)\) for all~\(\sit\in\sits\), and we collect these maps in the set~\(\supermartins\).
We call a non-negative supermartingale \(\test\) for~\(\frcstsystem\) such that \(\test(\init)=1\) a \emph{test supermartingale} for~\(\frcstsystem\), and we collect these in the set~\(\tests\).
\end{definition}

Readers familiar with the field of algorithmic randomness know that we mustn't allow Sabine to select just any allowable betting strategy---or test supermartingale.
Otherwise, the corresponding notion of randomness wouldn't make much sense because, for one thing, no path~\(\pth\in\pths\) would be random for the constant forecast~\(\nicefrac{1}{2}\).
This issue's typically resolved by restricting Sabine's betting strategies to a countable class of `effectively implementable' ones.
In \sectionref{sec:comp} we'll explain what `effectively implementable' means, but let's first have a look at how to devise a notion of randomness for an infinite sequence~\((\frcst_1,x_1,\dots,\frcst_n,x_n,\dots)\) of forecasts~\(\frcst_n\) and subsequent outcomes~\(x_n\)  when adopting a \emph{prequential} perspective.

\subsection{The Prequential Approach}
Again, we start by introducing a bit of notation.
An infinite sequence~\((\frcst_1,x_1,\dots,\frcst_n,x_n,\dots)\in\preqpths\) of rational forecasts~\(\frcst_n\) and subsequent outcomes~\(x_n\) is called a \emph{prequential path} and generically denoted by~\(\preqpth\).\footnote{We limit ourselves to rational forecasts in this prequential setting and draw attention to this restriction by using a subscript~\(r\); a rational forecasting system is for example denoted by~\(\ratfrcstsystem\), and the set of all rational interval forecasts by~\(\ratfrcsts\). In \sectionref{sec:randomness}, we'll provide some explanation and motivation for this restriction.}
An infinite sequence of rational forecasts~\((\frcst_1,\dots,\frcst_n,\dots)\in\frcstseqs\) is generically denoted by~\(\frcstseq\).
A finite sequence of rational forecasts and outcomes~\((\frcst_1,x_1,\dots,\frcst_n,x_n)\in\preqsits\), with \(\preqsits\coloneqq\bigcup_{k\in\naturalswithzero}\group{\smash{\ratfrcsts\times\posspace}}^k\), is called a \emph{prequential situation}, is generically denoted by~\(\preqsit\) and has length~\(\abs{\preqsit}=n\).
A finite sequence of rational forecasts~\((\frcst_1,\dots,\frcst_n)\in\frcstsits\), with \(\frcstsits\coloneqq\bigcup_{k\in\naturalswithzero}\ratfrcsts^k\), is generically denoted by~\(\frcstsit\) and has length~\(\abs{\frcstsit}=n\).
For any~\(k\in\naturalswithzero\), \(\preqpthto{k}=\group{I_1,x_1,\dots,I_k,x_k}\), and similarly for infinite sequences of rational forecasts~\(\frcstseq\in\frcstseqs\) and for prequential situations~\(\preqsit\in\preqsits\) with~\(k\leq\abs{\preqsit}\).
Furthermore, for any~\(k\in\naturalswithzero\), \(\frcstseqat{k}=I_k\), and similarly for a finite sequence of rational forecasts~\(\frcstsit\in\frcstsits\) with~\(k\leq\abs{\frcstsit}\).
The empty prequential situation~\(\preqsitto{0}=\group{}\) is denoted also by~\(\init\).

For ease of notation, we won't differentiate between~\(\preqpth\in\preqpths\) and \((\frcstseq,\pth)\in\frcstseqs\times\pths\).
In the same spirit, we won't differentiate between~\(\preqsit\in\preqsits\) and \((\frcstsit,\sit)\in\bigcup_{n\in\naturalswithzero}\ratfrcsts^n\times\posspace^n\).
The concatenation of a situation~\(\sit\in\sits\) and an outcome \(x\in\posspace\) is denoted by~\(\sit x\), the concatenation of a finite sequence of rational forecasts~\(\frcstsit\in\frcstsits\) and a rational forecast~\(\ratfrcst\in\ratfrcsts\) by~\(\frcstsit\ratfrcst\), and the concatenation of a prequential situation~\(\preqsit\in\preqsits\), a rational forecast~\(\ratfrcst\in\ratfrcsts\) and an outcome \(x\in\posspace\) by~\(\preqsit\ratfrcst x\).
In this way, for any~\(\preqsit=\preqsitrep=(\frcst_1,x_1,\dots,\frcst_n,x_n)\in\preqsits\), \(\ratfrcst\in\ratfrcsts\) and \(x\in\posspace\), we have that \(\preqsit\ratfrcst x=(\frcst_1,x_1,\dots,\frcst_n,x_n,\ratfrcst,x)=(\frcstsit\ratfrcst,\sit x)\in\preqsits\).

In the prequential setting, it's not assumed that Frank's forecasts are produced by some underlying forecasting system.
Instead, as is (re)presented in the protocol, he's allowed to produce forecasts on the fly, so there's no need for Frank to provide forecasts in all situations that could occur.
To test whether Frank is doing a good job, Sabine here too engages in a betting game, only now she has to define a strategy that specifies an allowed change in capital for all possible successions of rational forecasts (that could have been chosen by Frank) and outcomes (that could have been revealed by Reality), that is, she has to specify a possible change in capital for all prequential situations~\(\preqsit\in\preqsits\); she's again prohibited from borrowing money and assumed to start with unit capital.
Her prequential betting strategies are formalised as follows; as announced in the Introduction, we borrow the underlying idea as well as the terminology from Refs.~\cite{Dawid1999,Vovk2010}.

\begin{definition}
A real-valued map~\(\superfarth\colon\preqsits\to\reals\) is called a \emph{superfarthingale} if \allowbreak\(\uex_{\ratfrcst}(\superfarth(\preqsit\ratfrcst\andoutcome))\leq \superfarth(\preqsit)\) for all~\(\preqsit\in\preqsits\) and \(\ratfrcst\in\ratfrcsts\), and these maps are collected in the set~\(\superfarths\).
We call a non-negative superfarthingale~\(\superfarth\geq0\) such that \(\superfarth(\init)=1\) a \emph{test superfarthingale}.
\end{definition}

Here too, to obtain a sensible prequential notion of randomness, we need to restrict Sabine's betting strategies to a countable set, and we'll do so by specifying what it means for a betting strategy to be `effectively implementable'.

\section{Effective Objects}\label{sec:comp}
To define what it means for a mathematical object to be effectively implementable, we turn our attention to the field of computability theory.
As its basic objects, it considers natural maps~\(\phi\colon\naturals\to\naturals\).
Such a natural map~\(\phi\) is called \emph{recursive} if it can be computed by a Turing machine; this means that there's a Turing machine that, when given a natural number~\(n\in\naturals\), outputs the natural number~\(\phi(n)\in\naturals\).
By the Church-Turing thesis, this is equivalent to the existence of a finite algorithm that outputs \(\phi(n)\in\naturals\) in a finite number of steps when given \(n\in\naturals\) as an input.
Via encoding, this notion of effectiveness is extended to all rational maps~\(q\colon\mathcal{D}\to\rationals\), where \(\mathcal{D}\) denotes any countably infinite set whose elements can be encoded by the natural numbers; the choice of encoding isn't important, provided we can algorithmically decide whether a natural number is an encoding of an object and, if this is the case, we can find an encoding of the same object with respect to the other encoding \cite[p.~xvi]{Alexander2017book}.
By the Church-Turing thesis, a rational map~\(q\colon\mathcal{D}\to\rationals\) is then \emph{recursive} if there's some finite algorithm that outputs the rational number~\(q(d)\in\rationals\) in a finite number of steps, when it's given \(d\in\mathcal{D}\) as an input.
In line with the approach in Ref.~\cite{Pour-ElRichards2016}, we'll provide or describe an algorithm whenever we want to establish a map's recursive character.
In particular, since a finite number of algorithms can always be combined into one \cite{MichaelSipser2006}, a rational forecasting system~\(\ratfrcstsystem\in\ratfrcstsystems\) is called \emph{recursive} if there are two recursive maps~\(\underline{q},\overline{q}\colon\sits\to\rationals\) such that \(\lratfrcstsystem(\sit)=\underline{q}(\sit)\) and \(\uratfrcstsystem(\sit)=\overline{q}(\sit)\) for all~\(\sit\in\sits\).

Recursive maps can be used to provide notions of implementability for (extended) real-valued maps of the form \(r\colon\mathcal{D}\to\reals\cup\set{+\infty}\), whose co-domain isn't countably infinite.
Such a map~\(r\) is called \emph{lower semicomputable} if there's some recursive rational map~\(q\colon\smash{\mathcal{D}\times\naturalswithzero}\to\rationals\) such that \(q(d,n)\leq q(d,n+1)\) and \(\lim_{m\to\infty}q(d,m)=r(d)\) for all \(d\in\mathcal{D}\) and \(n\in\naturalswithzero\); as a shorthand notation, we'll then write \(q(d,\bolleke)\lsc r(d)\).
Any such map~\(q\) that witnesses the lower semicomputability of the map~\(r\) in the above sense, will also be called a \emph{code} for~\(r\).
We may always assume that this approximation from below is strictly increasing.

\begin{lemma} \label{lem:lsc}
An extended real map~\(r\colon\mathcal{D}\to\reals\cup\set{+\infty}\) is lower semicomputable if and only if there's some recursive rational map~\(q\colon\mathcal{D}\times\naturalswithzero\to\rationals\) such that \(\lim_{m\to\infty}q(d,m)=r(d)\) and \(q(d,n)<q(d,n+1)\) for all~\(d\in\mathcal{D}\) and \(n\in\naturalswithzero\).
\end{lemma}

\begin{proof}
The `if'-part is obvious.
For the `only if'-part, consider a recursive rational map \(q'\colon\mathcal{D}\times\naturalswithzero\to\rationals\) such that \(q'(d,\bolleke)\lsc r(d)\) for all~\(d\in\mathcal{D}\).
Define \(q\colon\mathcal{D}\times\naturalswithzero\to\rationals\) by~\(q(d,n)\coloneqq q'(d,n)-2^{-n}\) for all~\(d\in\mathcal{D}\) and \(n\in\naturalswithzero\).
Then \(\lim_{m\to\infty}q(d,m)=\lim_{m\to\infty}q'(d,m)=r(d)\) and \(q(d,n)<q'(d,n)-2^{-(n+1)}\leq q'(d,n+1)-2^{-(n+1)}=q(d,n+1)\) for all~\(d\in\mathcal{D}\) and \(n\in\naturalswithzero\).
\end{proof}
\noindent This also provides a proof for the following statement.

\begin{corollary} \label{cor:lsc}
There's a single algorithm that, upon input of a code for a lower semicomputable extended real map~\(r\colon\mathcal{D}\to\reals\cup\set{+\infty}\), outputs a recursive rational map \(q\colon\smash{\mathcal{D}\times\naturalswithzero}\to\rationals\) such that \(\lim_{m\to\infty}q(d,m)=r(d)\) and \(q(d,n)<q(d,n+1)\) for all~\(d\in\mathcal{D}\) and \(n\in\naturalswithzero\).
\end{corollary}

We'll also consider a stronger notion of effective implementability: a real map~\(r\colon\mathcal{D}\to\reals\) is called \emph{computable} if there's some recursive rational map~\(q\colon\mathcal{D}\times\naturalswithzero\to\rationals\) such that \(\abs{r(d)-q(d,n)}\leq2^{-n}\) for all~\(d\in\mathcal{D}\) and \(n\in\naturalswithzero\).
In particular, a forecasting system~\(\frcstsystem\in\frcstsystems\) is called computable if there are two recursive rational maps~\(\underline{q},\overline{q}\colon\sits\times\naturalswithzero\to\rationals\) such that \(\abs{\lfrcstsystem(\sit)-\underline{q}(\sit,n)}\leq2^{-n}\) and \(\abs{\ufrcstsystem(\sit)-\overline{q}(\sit,n)}\leq2^{-n}\) for all~\(\sit\in\sits\) and \(n\in\naturalswithzero\).

\section{Martin-Löf and Game-Randomness}\label{sec:randomness}

\subsection{The Standard Approach}
To get to a first notion of randomness, in the standard setting, we impose lower semicomputability on Sceptic Sabine's betting strategies---the test supermartingales---so as to obtain our `imprecise-probabilistic' martingale-theoretic version of Martin-Löf randomness \cite[Definition~2]{CoomanBock2021}.
We refer to our earlier work \cite{CoomanBock2021} for an extensive discussion of this type of randomness, its properties, and reasons for introducing it.

\begin{definition} \label{def:ml}
A path~\(\pth\in\pths\) is Martin-Löf random for a forecasting system~\(\frcstsystem\in\frcstsystems\) if \(\limsup_{n\to\infty}\test(\pthton)<\infty\) for all lower semicomputable test supermartingales \(\test\in\tests\).
\end{definition}

We can give a more prequential flavour to this randomness notion, but to do so, we require some more terminology.
With any infinite sequence of outcomes~\(\pth\in\pths\) and forecasting system~\(\frcstsystem\in\frcstsystems\), we associate the infinite sequence of forecasts~\(\frcstsystem\sqgroup{\pth}\coloneqq(\frcstsystem(\pthto{0}),\allowbreak\frcstsystem(\pthto{1}),\frcstsystem(\pthto{2}),\dots)\).
Similarly, we associate with any finite sequence of outcomes \(\sit\in\sits\) and forecasting system \(\frcstsystem\in\frcstsystems\) the finite sequence of forecasts~\(\frcstsystem\sqgroup{\sit}\coloneqq(\frcstsystem(\sitto{0}),\allowbreak\frcstsystem(\sitto{1}),\dots,\frcstsystem(\sitto{\abs{\sit}-1}))\).
This allows us to check the \emph{compatibility} of a forecasting system~\(\frcstsystem\in\frcstsystems\) with a given infinite sequence~\(\preqpth=\preqpthrep\in\realpreqpths\) of forecasts and outcomes, in the sense that \(\frcstsystem\) emits the same forecasts based on the observed outcomes~\(\pth\) in~\(\preqpth\) as the forecasts~\(\frcstseq\) that are present in~\(\preqpth\): we say that \(\frcstsystem\) is \emph{compatible} with~\(\preqpth\) if \(\frcstsystem\sqgroup{\pth}=\frcstseq\), that is, if \(\frcstsystem(\pthton)=\frcstseqat{n+1}\) for all~\(n\in\naturalswithzero\).
If the forecasting system~\(\frcstsystem\) produces more conservative forecasts along~\(\pth\) compared to \(\frcstseq\), that is, if \(\frcstseqatnplus\subseteq\frcstsystem(\pthton)\) for all~\(n\in\naturalswithzero\), then we say that~\(\frcstsystem\) is more \emph{conservative} (or less \emph{informative}) on~\(\preqpth=\preqpthrep\).
Similarly, we say that a forecasting system~\(\frcstsystem\) is \emph{compatible} with a prequential situation~\(\preqsit=\preqsitrep\in\preqsits\) if \(\frcstsystem(\sitto{n})=\frcstsitat{n+1}\) for all~\(0\leq n\leq\abs{\preqsit}-1\).
\definitionref{def:ml} can now be adapted to this new context as follows.

\begin{definition}\label{def:ml:alt}
We'll call a sequence~\(\preqpth=(\frcstseq,\pth)\in\realpreqpths\) of interval forecasts and outcomes \emph{Martin-Löf random} if there's some forecasting system~\(\frcstsystem\) that's compatible with~\(\preqpth\) such that \(\pth\) is Martin-Löf random for~\(\frcstsystem\).
\end{definition}

Before introducing an `imprecise-probabilistic' and martingale-theoretic prequential notion of randomness that's inspired by Vovk and Shen's work \cite{Vovk2010}, let's now first argue why we restrict our attention to \emph{rational} forecasts in the prequential setting.
First of all, compared to their approach in Ref.~\cite{Vovk2010}, it allows us to employ a technically less involved version of effective implementability that results in simpler proofs.
Secondly, and perhaps more importantly, we intend to compare our standard and prequential notions of randomness, and, as the following proposition shows, rational forecasts are enough to capture the essence of randomness in the standard setting.

\begin{proposition}\label{prop:rational:suffice}
For every non-degenerate computable forecasting system~\(\frcstsystem\in\frcstsystems\) there's a recursive rational forecasting system~\(\ratfrcstsystem\in\ratfrcstsystems\), with~\(\frcstsystem\subseteq\ratfrcstsystem\), such that a path~\(\pth\in\pths\) is Martin-Löf random for~\(\frcstsystem\) if and only if it's Martin-Löf random for~\(\ratfrcstsystem\).
\end{proposition}

\begin{proof}
Since \(\frcstsystem\) is computable, there are two recursive rational maps~\(\underline{q},\overline{q}\colon\sits\times\naturalswithzero\to\rationals\) such that
\begin{multline}\label{eq:comp:one}
\abs[\big]{\lfrcstsystem(\sit)-\underline{q}(\sit,n)}\leq2^{-n}
\text{ and }
\abs[\big]{\ufrcstsystem(\sit)-\overline{q}(\sit,n)}\leq2^{-n}
\text{ for all~\(\sit\in\sits\) and \(n\in\naturalswithzero\)}.
\end{multline}
By \lemmaref{lem:supermartin:bounded:above} in Appendix~\ref{sec:appendix}, since \(\frcstsystem\) is also assumed to be non-degenerate, we know there's a recursive natural map~\(C\colon\sits\to\naturals\) such that \(\test(\sit)\leq C(\sit)\) for all~\(\sit\in\sits\) and \(\test\in\tests\).
We'll now use these recursive maps to define an appropriate recursive rational approximation of \(\frcstsystem\).
As a first step, let \(N\colon\sits\to\naturalswithzero\) be defined as
\begin{equation*}
N(\sit)\coloneqq\min\cset[\bigg]{n\in\naturalswithzero}{2^{-n}\leq\frac{2^{-\abs{\sit}}}{\max\{C(\sit1),C(\sit0)\}+2}}
\end{equation*}
for all~\(\sit\in\sits\).
This map is recursive because \(C\) is.
Now, for any~\(\sit\in\sits\), let \(\ratfrcstsystem\in\ratfrcstsystems\) be defined by
\begin{align*}
\lratfrcstsystem(\sit)&\coloneqq\max\set[\big]{0,\underline{q}(\sit,N(\sit)+1)-2^{-(N(\sit)+1)}}
\shortintertext{and}
\uratfrcstsystem(\sit)&\coloneqq\min\set[\big]{1,\overline{q}(\sit,N(\sit)+1)+2^{-(N(\sit)+1)}}.
\end{align*}
By \equationref{eq:comp:one}, \(\underline{q}(\sit,N(\sit)+1)-2^{-(N(\sit)+1)}\leq\lfrcstsystem(\sit)\), and hence, since \(0\leq\lfrcstsystem(\sit)\), also \(\lratfrcstsystem(\sit)\leq\lfrcstsystem(\sit)\),  for all~\(\sit\in\sits\).
By \equationref{eq:comp:one}, it also holds for all~\(\sit\in\sits\) that \(\lfrcstsystem(\sit)\leq \underline{q}(\sit,N(\sit)+1)+2^{-(N(\sit)+1)}\), and therefore \(\lfrcstsystem(\sit)-2^{-N(\sit)}\leq\underline{q}(\sit,N(\sit)+1)-2^{-(N(\sit)+1)}\leq\lratfrcstsystem(\sit)\).
We conclude that
\begin{equation}\label{eq:ineq:one}
\lfrcstsystem(\sit)-2^{-N(\sit)}\leq\lratfrcstsystem(\sit)\leq\lfrcstsystem(\sit)\text{ for all } \sit\in\sits.
\end{equation}
In a similar fashion, we can show that
\begin{equation}\label{eq:ineq:two}
\ufrcstsystem(\sit)\leq\uratfrcstsystem(\sit)\leq\ufrcstsystem(\sit)+2^{-N(\sit)}
\text{ for all } \sit\in\sits.
\end{equation}
As a result, we already find that \(\frcstsystem\subseteq\ratfrcstsystem\).
Proposition~10 in Ref.~\cite{CoomanBock2021} then tells us that a path~\(\pth\in\pths\) is Martin-Löf random for~\(\frcstsystem\) only if it's Martin-Löf random for~\(\ratfrcstsystem\).
It remains to prove the `if'-direction, so assume that \(\pth\) is Martin-Löf random for~\(\ratfrcstsystem\) and assume towards contradiction that there's some lower semicomputable test supermartingale \(\test\in\tests\) for which \(\limsup_{n\to\infty}\test(\pthton)=\infty\).
Define the map~\(\test'\colon\sits\to\reals\) as
\begin{equation*}
\test'(\sit)
\coloneqq\frac{\test(\sit)+2^{-\abs{\sit}+1}}{3} \text{ for all } \sit\in\sits.
\end{equation*}
Clearly, \(\limsup_{n\to\infty}\test'(\pthton)=\frac{1}{3}\limsup_{n\to\infty}\test(\pthton)=\infty\), so we're done if we can prove that \(\test'\in\tests[\ratfrcstsystem]\) and that \(\test'\) is lower semicomputable.
\(\test'\) starts with unit capital since \(\frac{\test(\init)+2}{3}=1\), is non-negative since \(\test'(\sit)\geq\frac{\test(\sit)}{3}\geq0\) for all~\(\sit\in\sits\), and is lower semicomputable because \(\test\) is.
We complete the proof by proving its supermartingale character.
Fix any~\(\sit\in\sits\).
If \(\test'(\sit1)\geq\test'(\sit0)\), then
\begin{align*}
\uex_{\ratfrcstsystem(\sit)}(\test'(\sit\andoutcome))
\overset{\eqref{def:uex}}&{=}
\uratfrcstsystem(\sit)\test'(\sit1)+\group[\big]{1-\uratfrcstsystem(\sit)}\test'(\sit0) \\
\overset{\eqref{eq:ineq:two}}&{\leq}
\group[\big]{\ufrcstsystem(\sit)+2^{-N(\sit)}}\test'(\sit1)+\group[\big]{1-\ufrcstsystem(\sit)}\test'(\sit0) \\
\overset{\eqref{def:uex}}&{=}
\uexfrcstsystem(\test'(\sit\andoutcome))+2^{-N(\sit)}\test'(\sit1) \\
\overset{\text{\ref{axiom:coherence:homogeneity},\ref{axiom:coherence:constantadditivity}}}&{=}
\frac{\uexfrcstsystem(\test(\sit\andoutcome))+2^{-\abs{\sit}}}{3}+2^{-N(\sit)}\test'(\sit1) \\
&\leq\frac{\test(\sit)+2^{-\abs{\sit}}}{3}+\frac{2^{-\abs{\sit}}}{\max\{C(\sit1),C(\sit0)\}+2}\frac{\test(\sit1)+2}{3} \\
&\leq\frac{\test(\sit)+2^{-\abs{\sit}}}{3}+\frac{2^{-\abs{\sit}}}{3} \\
&=\frac{\test(\sit)+2^{-\abs{\sit}+1}}{3}
=\test'(\sit).
\end{align*}
Otherwise, if \(\test'(\sit1)<\test'(\sit0)\), then
\begin{align*}
\uex_{\ratfrcstsystem(\sit)}(\test'(\sit\andoutcome))
\overset{\eqref{def:uex}}&{=}
\lratfrcstsystem(\sit)\test'(\sit1)+\group[\big]{1-\lratfrcstsystem(\sit)}\test'(\sit0) \\
\overset{\eqref{eq:ineq:one}}&{\leq}
\lfrcstsystem(\sit)\test'(\sit1)+\group[\big]{1-\lfrcstsystem(\sit)+2^{-N(\sit)}}\test'(\sit0)\\
\overset{\eqref{def:uex}}&{=}
\uexfrcstsystem(\test'(\sit\andoutcome))+2^{-N(\sit)}\test'(\sit0) \\
\overset{\text{\ref{axiom:coherence:homogeneity},\ref{axiom:coherence:constantadditivity}}}&{=}
\frac{\uexfrcstsystem(\test(\sit\andoutcome))+2^{-\abs{\sit}}}{3}+2^{-N(\sit)}\test'(\sit0) \\
&\leq\frac{\test(\sit)+2^{-\abs{\sit}}}{3}+\frac{2^{-\abs{\sit}}}{\max\{C(\sit1),C(\sit0)\}+2}\frac{\test(\sit0)+2}{3} \\
&\leq\frac{\test(\sit)+2^{-\abs{\sit}}}{3}+\frac{2^{-\abs{\sit}}}{3} 
=\test'(\sit),
\end{align*}
so we're done.
\end{proof}

\subsection{A Prequential (Martingale-Theoretic) Approach}
To obtain a truly prequential \allowbreak imprecise-probabilistic martingale-theoretic notion of randomness, we mimic Vovk and Shen's approach \cite{Vovk2010}, and proceed by imposing lower semicomputability on Sabine's prequential betting strategies---which we've called test superfarthingales.
Contrary to their approach, we won't allow the test superfarthingales to be infinite-valued as a way to take care of conditional probability zero; instead, to deal with this issue, we explicitly restrict our attention to prequential paths \(\preqpth=\preqpthrep\in\preqpths\) that don't allow zero probability jumps, i.e., for which \(\frcstseqat{n}\neq1-\pthatn\) for all~\(n\in\naturals\), and which we'll call \emph{non-degenerate} prequential paths.
Analogously, we'll call a prequential situation~\(\preqsit=\preqsitrep\in\preqsits\) \emph{non-degenerate} if \(\frcstsitat{m}\neq1-\sitat{m}\) for all~\(1\leq m\leq\abs{\sit}\).

\begin{definition}
We call a sequence~\(\preqpth=\preqpthrep\in\preqpths\) of rational forecasts and outcomes \emph{game-random} if it's non-degenerate and if all lower semicomputable test {superfarthingales}~\(\superfarth\in\superfarths\) satisfy~\(\limsup_{n\to\infty}\superfarth(\preqpthton)<\infty\).
\end{definition}
\noindent In the remainder, we intend to explore how this new prequential randomness notion compares to our notion of Martin-Löf randomness.
We'll start by comparing definitions to uncover which (prequential) paths are(n't) random for both notions, and will then show that these definitions result in (almost) equivalent randomness notions when we restrict our attention to \emph{recursive} rational forecasting systems on the standard approach.
This endeavour can be seen as a continuation (and generalisation) of the discussion in Section~4 of Ref.~\cite{Vovk2010}, where Vovk and Shen prove that a standard and a prequential precise-probabilistic approach to randomness coincide for non-degenerate computable forecasting systems.
Afterwards, we'll compare a few basic properties for both imprecise-probabilistic notions, where we'll be especially concerned with whether (and which) computability restrictions are necessary for these properties to hold.

\section{Comparing Both Randomness Notions}\label{sec:comparison}

\subsection{Game-Randomness Implies Martin-Löf Randomness} \label{subsec:relation}
Any prequential path that's game-random is also Martin-Löf random.
Game-randomness is a therefore at least as \emph{strong} a randomness notion as Martin-Löf randomness.

\begin{proposition}\label{prop:rational:equivalencedirection:one}
Consider any infinite sequence of interval forecasts and outcomes~\(\preqpth=(\frcstseq,\pth)\in\preqpths\) that's game-random.
Then the infinite sequence of outcomes~\(\pth\) is Martin-Löf random for any rational forecasting system~\(\ratfrcstsystem\in\ratfrcstsystems\) that's compatible with~\(\preqpth\).
\end{proposition}

\begin{proof}
Consider any rational forecasting system~\(\ratfrcstsystem\in\ratfrcstsystems\) that's compatible with~\(\preqpth\) (which is non-degenerate by assumption) and assume towards contradiction that there's a lower semicomputable test supermartingale \(\test\in\tests[\ratfrcstsystem]\) such that \(\limsup_{n\to\infty}\test(\pthton)=\infty\); we can assume \(\test\) to be positive.
We'll now construct a lower semicomputable test superfarthingale~\(\superfarth'\in\superfarths\) in such a way that \(\superfarth'(\ratfrcstsystem[\sit],\sit)=\test(\sit)\) for all~\(\sit\in\sits\) for which \((\ratfrcstsystem[\sit],\sit)\) is non-degenerate, for which then of course \(\limsup_{n\to\infty}\superfarth'(\preqpthto{n})=\infty\).

Define the lower semicomputable map~\(\superfarth\colon\preqsits\to\reals\) by~\(\superfarth(\frcstsit,\sit)\coloneqq\test(\sit)\) for all~\((\frcstsit,\sit)\in\preqsits\).
By construction, \(\superfarth(\ratfrcstsystem[\bolleke],\bolleke)\colon\sits\to\reals\) is a positive test supermartingale for~\(\ratfrcstsystem\).
Invoking \lemmaref{lem:uniform:program} in Appendix~\ref{sec:appendix}, we then obtain a lower semicomputable test superfarthingale~\(\superfarth'\in\superfarths\) such that \(\superfarth'(\ratfrcstsystem[\sit],\sit)=\test(\sit)\) for all~\(\sit\in\sits\) for which \((\ratfrcstsystem[\sit],\sit)\) is non-degenerate.
\end{proof}

Conversely, any Martin-Löf random path~\(\pth\in\pths\) is also game-random, provided we impose recursiveness on the forecasting systems \(\ratfrcstsystem\in\ratfrcstsystems\) and non-degeneracy on the prequential paths \((\ratfrcstsystem\sqgroup{\pth},\pth)\in\preqpths\).

\begin{proposition}\label{prop:rational:equivalencedirection:two}
Consider any recursive rational forecasting system~\(\ratfrcstsystem\in\ratfrcstsystems\), and any path \(\pth\in\pths\).
If \(\pth\) is Martin-Löf random for~\(\ratfrcstsystem\) and \((\ratfrcstsystem\sqgroup{\pth},\pth)\) is non-degenerate, then the prequential path~\((\ratfrcstsystem\sqgroup{\pth},\pth)\) is game-random.
\end{proposition}

\begin{proof}
Since \((\ratfrcstsystem\sqgroup{\pth},\pth)\) is non-degenerate, assume towards contradiction that there's some lower semicomputable test superfarthingale~\(\superfarth\in\superfarths\) such that \(\limsup_{n\to\infty}\superfarth(\ratfrcstsystem\sqgroup{\pthton},\pthton)=\infty\).
Let \(\test\colon\sits\to\reals\) be defined as \(\test(\sit)\coloneqq\superfarth(\ratfrcstsystem\sqgroup{\sit},\sit)\) for all~\(\sit\in\sits\).
Then it holds that \(\limsup_{n\to\infty}\test(\pthton)=\infty\).
So we're done if we can show that \(\test\in\tests[\ratfrcstsystem]\) and that \(\test\) is lower semicomputable.
Obviously, it holds that \(\test(\init)=1\) and \(\test\geq0\) because \(\superfarth(\init)=1\) and \(\superfarth\geq0\).
Furthermore, for any~\(\sit\in\sits\), it follows from the superfarthingale condition that \(\uexratfrcstsystem(\test(\sit\andoutcome))=\uexratfrcstsystem(\superfarth(\ratfrcstsystem[\sit]\ratfrcstsystem(\sit),\sit\andoutcome))\leq\superfarth(\ratfrcstsystem[\sit],\sit)=\test(\sit)\), so we conclude that \(\test\in\tests[\ratfrcstsystem]\).
Since \(\superfarth\) is assumed to be lower semicomputable, there's some recursive rational map~\(q\colon\preqsits\times\naturalswithzero\to\rationals\) such that \(q(\preqsit,\bolleke)\lsc\superfarth(\preqsit)\) for all~\(\preqsit\in\preqsits\).
Let the rational map~\(q'\colon\sits\times\naturalswithzero\to\rationals\) be defined as \(q'(\sit,n)\coloneqq q((\ratfrcstsystem[\sit],\sit),n)\) for all~\(\sit\in\sits\) and \(n\in\naturalswithzero\).
This is a recursive map since \(\ratfrcstsystem\) is assumed to be a recursive rational forecasting system.
By construction, \(q'(\sit,\bolleke)=q((\ratfrcstsystem[\sit],\sit),\bolleke)\lsc\superfarth(\ratfrcstsystem[\sit],\sit)=\test(\sit)\) for all~\(\sit\in\sits\), and therefore \(\test\) is lower semicomputable.
\end{proof}

The following example shows that the recursiveness of the rational forecasting systems~\(\ratfrcstsystem\in\ratfrcstsystems\) in the previous proposition can't be dropped, so game-randomness is a \emph{strictly stronger} randomness notion than Martin-Löf randomness, since there's at least one prequential path~\(\preqpth\in\preqpths\) that's Martin-Löf random but not game-random.

\begin{example}\label{ex:rational:equivalencedirection2}\rm
By Corollary~20 in Ref.~\cite{CoomanBock2021}, there's at least one path~\(\pth\in\pths\) that's Martin-Löf random for the stationary forecasting system~\(\nicefrac{1}{2}\); this path is then necessarily non-recursive.
Consider the rational forecasting system~\(\ratfrcstsystem\in\ratfrcstsystems\) defined by
\begin{equation*}
\ratfrcstsystem(\sit)\coloneqq\begin{cases}
\sqgroup{0,\nicefrac{1}{2}} &\text{if }\pthat{\abs{\sit}+1}=1 \\
\sqgroup{\nicefrac{1}{2},1} &\text{if }\pthat{\abs{\sit}+1}=0
\end{cases}
\text{ for all } \sit\in\sits,
\end{equation*}
which is non-recursive since \(\pth\) is.
Since \(\set{\nicefrac{1}{2}}\subseteq\ratfrcstsystem\), it follows from Proposition \(10\) in Ref. \cite{CoomanBock2021} that \(\pth\) is also Martin-Löf random for~\(\ratfrcstsystem\).
Meanwhile, \((\ratfrcstsystem[\pth],\pth)\) isn't game-random.
To see this, define the test superfarthingale~\(\superfarth\in\superfarths\) recursively by~\(\superfarth(\init)\coloneqq1\) and
\begin{align*}
\superfarth(\preqsit\ratfrcst x)\coloneqq\begin{cases}
2\superfarth(\preqsit) &\text{if }\ratfrcst=\sqgroup{0,\nicefrac{1}{2}}\text{ and }x=1\\
2\superfarth(\preqsit) &\text{if }\ratfrcst=\sqgroup{\nicefrac{1}{2},1}\text{ and }x=0\\
0 &\text{otherwise}
\end{cases}
\text{ for all }\preqsit\in\preqsits\text{, }\ratfrcst\in\ratfrcsts\text{ and }x\in\posspace.
\end{align*}
\(\superfarth\) is clearly recursive, and \(\superfarth(\ratfrcstsystem[\pthton],\pthton)=2^n\) for all~\(n\in\naturalswithzero\).
Consequently, it holds that \(\limsup_{n\to\infty}\superfarth(\ratfrcstsystem[\pthton],\pthton)=\infty\), and therefore \((\ratfrcstsystem[\pth],\pth)\) isn't game-random.
\end{example}

By combining the last two propositions, we obtain conditions under which both randomness notions coincide; these conditions mimic the ones in Corollary~1 of Ref.~\cite{Vovk2010}, which are required to obtain a similar equivalence in Vovk and Shen's precise-probabilistic setting.

\begin{theorem}\label{the:rational:equivalence}
Consider any non-degenerate recursive rational forecasting system \(\ratfrcstsystem\in\ratfrcstsystems\).
Then any path~\(\pth\in\pths\) is Martin-Löf random for~\(\ratfrcstsystem\) if and only if the prequential path~\((\ratfrcstsystem\sqgroup{\pth},\pth)\in\preqpths\) is game-random.
\end{theorem}
\noindent This also shows that if a path~\(\pth\in\pths\) is Martin-Löf random for a non-degenerate recursive rational forecasting system~\(\ratfrcstsystem\in\ratfrcstsystems\), then only the forecasts~\(\ratfrcstsystem[\pth]\) that are produced along~\(\pth\) matter, since the path~\(\pth\) is also Martin-Löf random for any other non-degenerate recursive rational forecasting system~\(\ratfrcstsystem'\in\ratfrcstsystems\) such that \(\ratfrcstsystem'[\pth]=\ratfrcstsystem[\pth]\).
This result is in line with Dawid's Weak Prequential Principle \cite{Dawid1999}, which states that any criterion for assessing the `agreement' between Forecaster Frank and Reality should depend only on the actual observed sequences \(\frcstseq=(\frcst_1,\dots,\frcst_n,\dots)\in\frcstseqs\) and \(\pth=(x_1,\dots,x_n,\dots)\in\pths\), and not on the strategies (if any) which might have produced these, such as a recursive rational forecasting system~\(\ratfrcstsystem\in\ratfrcstsystems\) for  which \(\ratfrcstsystem[\pth]=\frcstseq\) .

\subsection{Properties}\label{subsec:properties}
As a first property, similarly as for Martin-Löf randomness \cite{DowneyHirschfeldt2010,MartinLof1966}, we mention (and prove) the existence of a so-called \emph{universal} test superfarthingale~\(\unifarth\in\superfarths\) that conclusively tests the game-randomness of any non-degenerate prequential path~\(\preqpth\in\preqpths\).

\begin{theorem}\label{the:rational:universal}
There's a so-called \emph{universal} superfarthingale~\(\unifarth\) with the property that any non-degenerate prequential path~\(\preqpth\in\preqpths\) is game-random iff \(\limsup_{n\to\infty}\unifarth(\preqpthton)<\infty\).
\end{theorem}

\begin{proof}
Lemma~13 in Ref.~\cite{Vovk2010} states that there's a uniformly lower semicomputable sequence of maps~\(f_n\colon\preqsits\to\sqgroup{0,+\infty}\) that contains every lower semicomputable map \(f\colon\preqsits\to\sqgroup{0,+\infty}\).
The sequence~\(\group{f_n}_{n\in\naturalswithzero}\) contains all lower semicomputable positive test superfarthingales \(\superfarth\in\superfarths\), so it follows from \corollaryref{cor:uniform:program} in Appendix~\ref{sec:appendix} that there's a uniformly lower semicomputable sequence of test superfarthingales \(\superfarth_n\in\superfarths\) such that for every positive test superfarthingale~\(\superfarth'\in\superfarths\) there's some \(N\in\naturalswithzero\) such that
\begin{align*}
\superfarth_N(\preqsit)=\begin{cases}
\superfarth'(\preqsit) &\text{if }\preqsit\text{ is non-degenerate} \\
0 &\text{if }\preqsit\text{ is degenerate}
\end{cases}
\text{ for all }\preqsit\in\preqsits.
\end{align*}

Let \(\unifarth\colon\preqsits\to\reals\) be defined by~\(\unifarth(\preqsit)\coloneqq\sum_{n=0}^{\infty}2^{-n-1}\superfarth_n(\preqsit)\) for all~\(\preqsit\in\preqsits\).
Since \(\superfarth_n\geq0\) and \(\superfarth_n(\init)=1\) for all~\(n\in\naturalswithzero\), it follows that \(\unifarth\) is well-defined (although possibly infinite), \(\unifarth\geq0\) and \(\unifarth(\init)=1\).
To check that \(\unifarth\) is indeed real-valued, fix any prequential situation~\(\preqsit=\preqsitrep\in\preqsits\).
If \(\preqsit\) is degenerate, then \(\unifarth(\preqsit)=0\) because \(\superfarth_n(\preqsit)=0\) for all~\(n\in\naturalswithzero\) by \corollaryref{cor:uniform:program} in Appendix~\ref{sec:appendix}.
If \(\preqsit\) is non-degenerate, then we infer from \lemmaref{lem:superfarth:bound} in Appendix~\ref{sec:appendix} that there's some real number~\(B\in\reals\) such that \(\superfarth_n(\preqsit)\leq B\) for all~\(n\in\naturalswithzero\), and therefore \(\unifarth(\preqsit)\leq\sum_{n=0}^\infty2^{-n-1}B=B\).
A standard argument we won't repeat here shows that \(\unifarth\) is lower semicomputable as an infinite sum of uniformly lower semicomputable non-negative maps~\(\superfarth_n\).
To show that \(\unifarth\) is a superfarthingale, fix any~\(\preqsit\in\preqsits\) and any~\(\ratfrcst\in\ratfrcsts\), and observe that
\begin{align*}
\uex_{\ratfrcst}(\unifarth(\preqsit\ratfrcst\andoutcome))
&=\lim_{k\to\infty}\uex_{\ratfrcst}\group[\bigg]{\sum_{n=0}^k2^{-n-1}\superfarth_n(\preqsit\ratfrcst\andoutcome)} \\
\overset{\text{\ref{axiom:coherence:homogeneity},\ref{axiom:coherence:subsupadditivity}}}&{\leq}
\sum_{n=0}^{\infty}2^{-n-1}\uex_{\ratfrcst}(\superfarth_n(\preqsit\ratfrcst\andoutcome)) \\
&\leq\sum_{n=0}^{\infty}2^{-n-1}\superfarth_n(\preqsit)
=\unifarth(\preqsit),
\end{align*}
where the first equality follows from \ref{axiom:coherence:uniform:convergence}, the real-valuedness of \(\unifarth\) and the non-negativity of \(\superfarth_n\) for all~\(n\in\naturalswithzero\), and the second inequality follows from the superfarthingale property for all~\(\superfarth_n\).
We conclude that \(\unifarth\) is a lower semicomputable test superfarthingale. 

We claim that \(\unifarth\) is a universal superfarthingale in the sense of the theorem.
Consider any non-degenerate prequential path~\(\preqpth\in\preqpths\).
The `only if'-part is obvious: if \(\preqpth\) is game-random, then \(\limsup_{n\to\infty}\superfarth'(\preqpthton)<\infty\) for all lower semicomputable test superfarthingales \(\superfarth'\in\superfarths\), and therefore also for~\(\unifarth\).
For the `if'-part, assume towards contradiction that there's some lower semicomputable test superfarthingale~\(\superfarth'\in\superfarths\) such that \(\limsup_{n\to\infty}\superfarth'(\preqpthton)=\infty\); we can assume \(\superfarth'\) to be positive: if it isn't, replace it with~\(\nicefrac{(\superfarth'+1)}{2}\).
We then know there's some \(N\in\naturalswithzero\) such that
\begin{equation*}
\superfarth_N(\preqsit)
=\begin{cases}
\superfarth'(\preqsit) &\text{if }\preqsit\text{ is non-degenerate} \\
0 &\text{if }\preqsit\text{ is degenerate}
\end{cases}
\text{ for all~\(\preqsit\in\preqsits\)}.
\end{equation*}
Hence, by the non-negativity of the test superfarthingales \(\superfarth_n\in\superfarths\) and the non-degeneracy of \(\preqpth\),
\begin{align*}
\limsup_{n\to\infty}\unifarth(\preqpthton)
&=\limsup_{n\to\infty}\sum_{k=0}^{\infty}2^{-k-1}\superfarth_k(\preqpthton) \\
&\geq\limsup_{n\to\infty}2^{-N-1}\superfarth_N(\preqpthton) \\
&=2^{-N-1}\limsup_{n\to\infty}\superfarth'(\preqpthton)
=\infty,
\end{align*}
so we're done.
\end{proof}

For Martin-Löf randomness, where the emphasis lies on the compatibility between a path and a forecasting system, we have that, for every forecasting system~\(\frcstsystem\in\frcstsystems\), there's at least one path~\(\pth\in\pths\) that's Martin-Löf random for~\(\frcstsystem\) \cite[Corollary 20]{CoomanBock2021}.
In the prequential setting, we have an analogous result for sequences of rational forecasts~\(\frcstseq\in\frcstseqs\) and sequences of outcomes~\(\pth\in\pths\).

\begin{proposition}\label{prop:rational:atleastone}
For every infinite sequence of rational interval forecasts~\(\frcstseq\in\frcstseqs\) there's at least one path~\(\pth\in\pths\) such that~\((\frcstseq,\pth)\in\preqpths\) is game-random.
\end{proposition}

\begin{proof}
Consider the universal superfarthingale~\(\unifarth\) from \theoremref{the:rational:universal}.
Assume that the path~\(\pth\) has been defined up to \(n\geq0\) entries such that \(1=\unifarth(\init)\geq\unifarth(\frcstseqto{1},\pthto{1})\geq\dots\geq\unifarth(\frcstseqton,\pthton)\) and \(\frcstseqat{m}\neq1-\pthat{m}\) for all~\(1\leq m\leq n\).
If \(\frcstseqat{n+1}=0\), let \(\pthatnplus\coloneqq0\).
Else, if \(\frcstseqat{n+1}=1\), let \(\pthatnplus\coloneqq1\).
In both cases, it holds by the superfarthingale property that
\begin{align*}
\unifarth(\frcstseqton,\pthton)&\geq\uex_{\frcstseqatnplus}(\unifarth(\frcstseqtonplus,\pthton\andoutcome)) \\
&=\frcstseqat{n+1}\unifarth(\frcstseqtonplus,\pthton1)+(1-\frcstseqat{n+1})\unifarth(\frcstseqtonplus,\pthton0)=\unifarth(\frcstseqtonplus,\pthton\frcstseqatnplus).
\end{align*}
Otherwise, if \(\frcstseqat{n+1}\notin\outcomes\), by the superfarthingale property and \ref{axiom:coherence:bounds}, we know that there's always some \(x\in\posspace\) such that
\[
\unifarth(\frcstseqton,\pthton)
\geq\uex_{\frcstseqatnplus}(\unifarth(\frcstseqtonplus,\pthton\andoutcome))
\geq\min\unifarth(\frcstseqtonplus,\pthton\andoutcome)
=\unifarth(\frcstseqtonplus,\pthton x),
\]
and then let \(\pthatnplus\coloneqq x\).
Invoking the axiom of dependent choice, we obtain a non-degenerate prequential path~\(\preqpth=\preqpthrep\in\preqpths\) such that \(\limsup_{n\to\infty}\unifarth(\preqpthton)\leq1\).
\end{proof}

In the next proposition and theorem, the required computability conditions on sequences of rational interval forecasts (in the prequential setting) differ from the ones on forecasting systems that are needed to obtain similar results in the standard setting \cite{CoomanBock2021}.
For example, any path~\(\pth\in\pths\) that's Martin-Löf random for a forecasting system~\(\frcstsystem\in\frcstsystems\) is also Martin-Löf random for any other more conservative forecasting system \cite[Proposition 10]{CoomanBock2021}.
Meanwhile, for a similar result to hold in the prequential setting, we need to restrict our attention to sequences of rational forecasts that are not only more conservative, but that also have a compatible recursive rational forecasting system.

\begin{proposition}\label{prop:rational:monotonicity}
Consider any recursive rational forecasting system~\(\ratfrcstsystem\in\ratfrcstsystems\) and any game-random prequential path~\(\preqpth=\preqpthrep\in\preqpths\).
If \(\ratfrcstsystem\) is more conservative on~\(\preqpth\), then \((\ratfrcstsystem\sqgroup{\pth},\pth)\) is game-random as well.
\end{proposition}

\begin{proof}
Consider any recursive rational forecasting system~\(\ratfrcstsystem\in\ratfrcstsystems\) that's more conservative on \(\preqpth\).
We can always consider a rational forecasting system~\(\ratfrcstsystem'\in\ratfrcstsystems\) that's compatible with~\(\preqpth\) such that \(\ratfrcstsystem'\subseteq\ratfrcstsystem\); note that we don't require computability here and that \(\preqpth=(\ratfrcstsystem'\sqgroup{\pth},\pth)\) is non-degenerate since \(\preqpth\) is game-random.
By \propositionref{prop:rational:equivalencedirection:one}, we then know that \(\pth\) is Martin-Löf random for~\(\ratfrcstsystem'\).
Consequently, by Proposition~10 in Ref.~\cite{CoomanBock2021}, since \(\ratfrcstsystem'\subseteq\ratfrcstsystem\), \(\pth\) is also Martin-Löf random for~\(\ratfrcstsystem\).
Since \(\preqpth=(\ratfrcstsystem'\sqgroup{\pth},\pth)\) is non-degenerate and \(\ratfrcstsystem'\subseteq\ratfrcstsystem\), the prequential path~\((\ratfrcstsystem\sqgroup{\pth},\pth)\) is non-degenerate, and hence, by \propositionref{prop:rational:equivalencedirection:two}, \((\ratfrcstsystem\sqgroup{\pth},\pth)\) is game-random too.
\end{proof}

The computability requirement in \propositionref{prop:rational:monotonicity} is not only sufficient, but also necessary.
This follows almost immediately from \exampleref{ex:rational:equivalencedirection2}.
It only remains to note that the stationary forecasting system~\(\nicefrac{1}{2}\) is rational and recursive, and hence, by \propositionref{prop:rational:equivalencedirection:two}, \((\nicefrac{1}{2},\pth)\) is game-random, while for the more conservative but non-recursive forecasting system~\(\ratfrcstsystem\), the prequential path~\((\ratfrcstsystem[\pth],\pth)\) isn't.

There are also prequential properties where the required computability conditions on the forecasts are less, rather than more, stringent.
If we restrict our attention for example to computable forecasting systems \(\frcstsystem\in\frcstsystems\), then the frequency of ones along a Martin-Löf random path, and along all so-called computably selected infinite subsequences, is bounded by the computable forecasting system \cite[Theorem 23]{CoomanBock2021}.
In the prequential setting, we have a similar result, but without any computability requirement on the infinite sequence of interval forecasts; in spirit, this result also generalises Dawid's ideas on calibration in Ref.~\cite{Dawid1982}.

\begin{theorem}\label{the:rational:freq}
Consider any infinite sequence of rational interval forecasts and outcomes~\(\preqpth\allowbreak=(\frcstseq,\pth)\in\preqpths\) and any recursive \emph{selection function}~\(\selection\colon\preqsits\times\ratfrcsts\to\outcomes\) such that~\(\sum_{k=0}^{\infty}\selection(\preqpthtok,\frcstseqatkplus)=\infty\).
If \(\preqpth\) is game-random, then
\begin{align*}
\liminf_{n\to\infty}\frac{\sum_{k=0}^{n-1}\selection(\preqpthtok,\frcstseqatkplus)\sqgroup{\pthatkplus-\min \frcstseqatkplus}}{\sum_{k=0}^{n-1}\selection(\preqpthtok,\frcstseqatkplus)}&\geq0
\shortintertext{and}
\limsup_{n\to\infty}\frac{\sum_{k=0}^{n-1}\selection(\preqpthtok,\frcstseqatkplus)\sqgroup{\pthatkplus-\max \frcstseqatkplus}}{\sum_{k=0}^{n-1}\selection(\preqpthtok,\frcstseqatkplus)}&\leq0.
\end{align*}
\end{theorem}

\begin{proof}
We'll give a proof for the first inequality, the proof for the second one is similar.
Assume towards contradiction that there's some real number~\(\epsilon\), with~\(0<\epsilon<1\), such that
\begin{equation*}
\liminf_{n\to\infty}\frac{\sum_{k=0}^{n-1}\selection(\preqpthtok,\frcstseqatkplus)\sqgroup{\pthatkplus-\min \frcstseqatkplus}}{\sum_{k=0}^{n-1}\selection(\preqpthtok,\frcstseqatkplus)}<-\epsilon.
\end{equation*}
Let the map~\(\superfarth\coloneqq\preqsits\to\reals\) be defined by
\begin{multline*}
\superfarth(\preqsit)
\coloneqq\prod_{k=0}^{\abs{\preqsit}-1}\sqgroup[\big]{1-\frac{\epsilon}{3}\selection(\preqsitto{k},\frcstsitat{k+1})[\sitat{k+1}-\min\frcstsitat{k+1}]}
\text{ for all~\(\preqsit=(\frcstsit,\sit)\in\preqsits\).}
\end{multline*}
We'll now show in a number of steps that \(\superfarth\) is a lower semicomputable test superfarthingale for which \(\limsup_{n\to\infty}\superfarth(\preqpthton)=\infty\), implying that \(\preqpth\) can't be game-random.

Trivially, \(\superfarth(\init)=1\), and also \(\superfarth\geq0\), because \(\epsilon<1\), \(\abs{\selection}\leq1\) and \(\abs{x-\min \ratfrcst}\leq1\) for all~\(x\in\posspace\) and \(\ratfrcst\in\ratfrcsts\).
Moreover, for any~\(\preqsit\in\preqsits\) and \(\ratfrcst\in\ratfrcsts\), we have that
\begin{align*}
\uex_{\ratfrcst}(\superfarth(\preqsit\ratfrcst\andoutcome))
\overset{\text{\ref{axiom:coherence:homogeneity}}}&{=}
\superfarth(\preqsit)\uex_{\ratfrcst}\group[\Big]{1+\frac{\epsilon}{3}\selection(\preqsit,\ratfrcst)[\min\ratfrcst-X]} \\
\overset{\text{\ref{axiom:coherence:homogeneity},\ref{axiom:coherence:constantadditivity}}}&{=}
\superfarth(\preqsit)\sqgroup[\Big]{1+\frac{\epsilon}{3}\selection(\preqsit,\ratfrcst)\uex_{\ratfrcst}(\min\ratfrcst-X)} \\
\overset{\text{\ref{axiom:coherence:constantadditivity}}}&{=}
\superfarth(\preqsit)\sqgroup[\Big]{1+\frac{\epsilon}{3}\selection(\preqsit,\ratfrcst)(\min\ratfrcst+\uex_{\ratfrcst}(-X))} \\
\overset{\text{\eqref{def:lex}}}&{=}
\superfarth(\preqsit),
\end{align*}
so we find that \(\superfarth\) is a test superfarthingale.
From the recursiveness of \(\selection\) and the rational-valuedness of the forecasts~\(\ratfrcst\in\ratfrcsts\) and outcomes~\(x\in\posspace\) it follows that \(\superfarth\) is recursive, and therefore lower semicomputable as well.
We conclude that \(\superfarth\) is a lower semicomputable test superfarthingale.

By assumption, for any~\(m,M\in\naturalswithzero\), there's some \(N>m\) such that \(\sum_{k=0}^{N-1}\selection(\preqpthto{k},\frcstseqat{k+1})\geq M\) and
\begin{equation}\label{eq:the:rational:freq}
\frac{\sum_{k=0}^{N-1}\selection(\preqpthtok,\frcstseqatkplus)\sqgroup{\pthatkplus-\min \frcstseqatkplus}}{\sum_{k=0}^{N-1}\selection(\preqpthtok,\frcstseqatkplus)}<-\epsilon.
\end{equation}
This will allow us to obtain a lower bound for~\(\superfarth(\preqpthto{N})\).
Since \(1-\frac{\epsilon}{3}\selection(\preqsit,\ratfrcst)\sqgroup{x-\min\ratfrcst}>\nicefrac{1}{2}\) for all~\(\preqsit\in\preqsits\), \(\ratfrcst\in\ratfrcsts\) and \(x\in\posspace\), it holds that \(\superfarth(\preqpthto{N})=\exp(K)\), with
\begin{equation*}
K
\coloneqq\sum_{k=0}^{N-1}\ln\group[\Big]{1-\frac{\epsilon}{3}\selection(\preqpthto{k},\frcstseqat{k+1})[\pthat{k+1}-\min\frcstseqat{k+1}]}.
\end{equation*}
Since \(\ln(1+x)\geq x-x^2\) for all~\(x>\nicefrac{-1}{2}\), we infer that
\begin{align*}
K
&\geq-\frac{\epsilon}{3}\sum_{k=0}^{N-1}\selection(\preqpthto{k},\frcstseqat{k+1})[\pthat{k+1}-\min\frcstseqat{k+1}]-\frac{\epsilon^2}{9}\sum_{k=0}^{N-1}\selection(\preqpthto{k},\frcstseqat{k+1})^2[\pthat{k+1}-\min\frcstseqat{k+1}]^2\\
\shortintertext{and, also taking into account \equationref{eq:the:rational:freq}, \(S^2=S\) and \([\pthat{k+1}-\min\frcstseqat{k+1}]^2\leq1\),}
&\geq\frac{\epsilon^2}{3}\sum_{k=0}^{N-1}\selection(\preqpthto{k},\frcstseqat{k+1})-\frac{\epsilon^2}{9}\sum_{k=0}^{N-1}\selection(\preqpthto{k},\frcstseqat{k+1}) \\
&=\frac{2\epsilon^2}{9}\sum_{k=0}^{N-1}\selection(\preqpthto{k},\frcstseqat{k+1}).
\end{align*}
Hence,
\begin{equation*}
\superfarth(\preqpthto{N})
\geq\exp\group[\bigg]{\frac{2\epsilon^2}{9}\sum_{k=0}^{N-1}\selection(\preqpthto{k},\frcstseqat{k+1})}
\geq\exp\group[\bigg]{\frac{2\epsilon^2}{9}M}.
\end{equation*}
After recalling that the inequality above holds for any~\(M\in\naturalswithzero\) and for arbitrarily large well-chosen \(N\in\naturalswithzero\), we conclude that \(\limsup_{n\to\infty}\superfarth(\preqpthto{n})=\infty\).
\end{proof}

\section{Conclusions and Future Work}
We've introduced an imprecise-probabilistic prequential notion of randomness, argued why we restrict our attention here to rational interval forecasts, and proved several properties of this randomness notion.
We're especially satisfied with having achieved equipping our \emph{standard} imprecise-probabilistic version of Martin-Löf randomness with a \emph{prequential} interpretation.

In future work, we intend to come closer to Vovk and Shen's work \cite{Vovk2010}, by allowing for real interval forecasts.
We'll try to achieve this by adopting a more involved notion of lower semicomputability that allows for real maps~\(r\colon\mathcal{D'}\to\reals\) whose domain \(\mathcal{D'}\) can be uncountable, such as the set~\(\group{\frcsts\times\posspace}^{\naturals}\).
We suspect that, in this continuous setting, the necessary conditions to obtain analogous results to the propositions and theorems in \sectionref{subsec:properties} will be different; for one thing, we expect the computability condition on the forecasting systems in \propositionref{prop:rational:monotonicity} to drop, which would arguably yield a more natural monotonicity property.

In line with Refs.~\cite{MartinLof1966,Schnorr1971,Schnorr1971book,Vovk2010}, we intend to explore whether we can equip a prequential imprecise-probabilistic (martingale-theoretic) randomness notion with a measure-theoretic characterisation.
We have already succeeded in doing so in the context of this paper, but decided to omit these results because of page limitations.
The answer to this question remains open for the more general prequential randomness notion that's alluded to in the previous paragraph.

Lastly, we wonder whether we can give a precise-probabilistic interpretation to our prequential imprecise-probabilistic notion of randomness. 
In the standard setting, we've shown \cite{Persiau2022} that a path~\(\pth\in\pths\) is Martin-Löf random for an interval forecast~\(\frcst\in\frcsts\) if and only if it's random for some precise forecasting system~\(\frcstsystem_p\in\frcstsystems\) that's compatible with~\(\frcst\), in the sense that \(\frcstsystem_p(\sit)\in\frcst\) for all~\(\sit\in\sits\).
In this prequential context, we might be able to interpret an infinite sequence of forecasts~\(\frcstseq=(\frcst_1,\dots,\frcst_n,\dots)\in\frcstseqs\) as bounds on precise forecasts, and say that a prequential path~\(\preqpthrep\in\preqpths\) is game-random if and only if \((p_1,\pthat{1},p_2,\pthat{2},\dots)\in\preqpths\) is game-random for some infinite sequence of probabilities~\((p_1,\dots,p_n,\dots)\in\frcstseqs\) such that \(p_i\in\frcst_i\) for all~\(i\in\naturals\).

\section*{Acknowledgments}
Work on this paper was supported by the Research Foundation -- Flanders (FWO), project numbers 11H5521N (for Floris Persiau) and 3G028919 (for Gert de Cooman).
Gert de Cooman's research was also partially supported by a sabbatical grant from Ghent University, and from the FWO, reference number K801523N.
He also wishes to express his sincere gratitude to Jason Konek, whose ERC Starting Grant ``Epistemic Utility for Imprecise Probability'' under the European Union’s Horizon 2020 research and innovation programme (grant agreement no.~852677) allowed him to make a sabbatical stay at Bristol University's Department of Philosophy, and to Teddy Seidenfeld, whose funding helped realise a sabbatical stay at Carnegie Mellon University's Department of Philosophy.

\section*{Author Contributions}
As alluded to in the acknowledgments section of Ref.~\cite{CoomanBock2021}, Gert has long had the idea of developing a prequential randomness notion that accommodates interval forecasts, and made some initial attempts at doing so before deciding to focus on introducing imprecision for the more standard randomness notion.
After a period of dormancy, the idea was taken up again by Floris, who researched and developed the results in this paper and wrote a first draft.
Gert then checked the results and made various suggestions, and both authors collaborated intensively on revising the initial version, which led to the present paper.

\appendix
\section{Technical Lemmas and Proofs}\label{sec:appendix}

\begin{lemma}\label{lem:superfarth:bound}
For any non-degenerate prequential situation~\(\preqsit=\preqsitrep\in\preqsits\) and any non-negative superfarthingale~\(\superfarth\in\superfarths\), \(\superfarth(\preqsit)\leq\prod_{k=1}^{\abs{\sit}}\frac{1}{\group{\max\frcstseqat{k}}^{\sitat{k}}\group{1-\min\frcstseqat{k}}^{1-\sitat{k}}}\superfarth(\init)\).
\end{lemma}

\begin{proof}
Consider any non-degenerate prequential situation~\(\preqsit\ratfrcst x\in\preqsits\).
If \(x=1\) then \(0<\max\ratfrcst\leq1\), and
\begin{align*}
\superfarth(\preqsit\ratfrcst x)
&\leq\frac{1}{\max\ratfrcst}\sqgroup[\big]{\max\ratfrcst\superfarth(\preqsit\ratfrcst1)+(1-\max\ratfrcst)\superfarth(\preqsit\ratfrcst0)} \\
&\leq\frac{1}{\max\ratfrcst}\uex_{\ratfrcst}(\superfarth(\preqsit\ratfrcst\andoutcome))
\leq\frac{1}{\max\ratfrcst}\superfarth(\preqsit).
\end{align*}
If \(x=0\) then \(0\leq\min\ratfrcst<1\), and
\begin{align*}
\superfarth(\preqsit\ratfrcst x)
&\leq\frac{1}{1-\min\ratfrcst}\sqgroup[\big]{\min\ratfrcst\superfarth(\preqsit\ratfrcst1)+(1-\min\ratfrcst)\superfarth(\preqsit\ratfrcst0)} \\
&\leq\frac{1}{1-\min\ratfrcst}\uex_{\ratfrcst}(\superfarth(\preqsit\ratfrcst\andoutcome))
\leq\frac{1}{1-\min\ratfrcst}\superfarth(\preqsit).
\end{align*}
Above, the first, second and third inequalities follow from the non-negativity of \(\superfarth\), \equationref{def:uex} and the superfarthingale property, respectively.
Hence,
\begin{equation*}
\superfarth(\preqsit\ratfrcst x)\leq\frac{1}{\group{\max\ratfrcst}^{x}\group{1-\min\ratfrcst}^{1-x}}\superfarth(\preqsit).
\end{equation*}
A simple induction argument now leads to the desired result.
\end{proof}

\begin{lemma}\label{lem:supermartin:bounded:above}
For every non-degenerate computable forecasting system~\(\frcstsystem\in\frcstsystems\) there's a recursive natural map~\(C\colon\sits\to\naturals\) such that for every test supermartingale \(\test\in\tests\) it holds that \(\test(\sit)\leq C(\sit)\) for all~\(\sit\in\sits\).
\end{lemma}

\begin{proof}
Define the map~\(C'\colon\sits\to\reals\) by letting \(C'(\sit)\coloneqq\prod_{k=1}^{\abs{\sit}}\frac{1}{\ufrcstsystem(\sitto{k-1})^{\sitat{k}}\group{1-\lfrcstsystem(\sitto{k-1})}^{1-\sitat{k}}}\) for all~\(\sit\in\sits\).
This map is real-valued, since \(0<\ufrcstsystem\) and \(\lfrcstsystem<1\) by the non-degeneracy of \(\frcstsystem\).
Since \(\frcstsystem\) is computable, \(C'\) is computable as well.
Let's now prove that \(\test(\sit)\leq C'(\sit)\) for all~\(\sit\in\sits\).
Fix any situation~\(\sit\in\sits\) and any~\(x\in\posspace\).
If \(x=1\), then
\begin{align*}
\test(\sit x)
&\leq\frac{1}{\ufrcstsystem(\sit)}\sqgroup[\big]{\ufrcstsystem(\sit)\test(\sit1)+(1-\ufrcstsystem(\sit))\test(\sit0)} \\
&\leq\frac{1}{\ufrcstsystem(\sit)}\uex_{\frcstsystem(\sit)}(\test(\sit\andoutcome))
\leq\frac{1}{\ufrcstsystem(\sit)}\test(\sit).
\end{align*}
If \(x=0\), then
\begin{align*}
\test(\sit x)
&\leq\frac{1}{1-\lfrcstsystem(\sit)}\sqgroup[\big]{\lfrcstsystem(\sit)\test(\sit1)+(1-\lfrcstsystem(\sit))\test(\sit0)} \\
&\leq\frac{1}{1-\lfrcstsystem(\sit)}\uex_{\frcstsystem(\sit)}(\test(\sit\andoutcome))
\leq\frac{1}{1-\lfrcstsystem(\sit)}\test(\sit).
\end{align*}
Above, the first, second and third inequalities follow from the non-negativity of \(\test\), \equationref{def:uex} and the supermartingale property, respectively.
Hence,
\[
\test(\sit x)\leq\frac{1}{\ufrcstsystem(\sit)^{x}\group[\big]{1-\lfrcstsystem(\sit)}^{1-x}}\test(\sit).
\]
A simple induction argument now shows that indeed \(\test(\sit)\leq C'(\sit)\) for all~\(\sit\in\sits\).

Since \(C'\) is a computable real map, there's a recursive rational map~\(q\colon\sits\times\naturalswithzero\to\rationals\) such that \(\abs{C'(\sit)-q(\sit,n)}\leq2^{-n}\) for all~\(\sit\in\sits\) and \(n\in\naturalswithzero\).
Let \(C\colon\sits\to\naturals\) be defined as \(C(\sit)\coloneqq\max\{1,\ceil{q(\sit,1)+1}\}\) for all~\(\sit\in\sits\), with~\(\ceil{\bolleke}\colon\reals\to\integers\) the ceiling function and \(\integers\) the set of integer numbers.
It's easy to see that \(C\) is natural-valued, positive and recursive.
Furthermore, we have that \(\test(\sit)\leq C'(\sit)\leq q(\sit,1)+\nicefrac{1}{2}\leq C(\sit)\) for all~\(\sit\in\sits\).
\end{proof}

\begin{lemma}\label{lem:uniform:program}
There's a single algorithm that, upon the input of a code for a lower semicomputable map~\(\superfarth\colon\preqsits\to\sqgroup{0,+\infty}\), outputs a code for a lower semicomputable test superfarthingale~\(\superfarth'\in\superfarths\) such that
\begin{enumerate}[label=\upshape(\roman*),leftmargin=*,itemsep=0pt]
\item \(\superfarth'(\preqsit)=0\) for all degenerate prequential situations~\(\preqsit\in\preqsits\); \label{lem:prop:one}
\item for any rational forecasting system~\(\ratfrcstsystem\in\ratfrcstsystems\), \(\superfarth'(\ratfrcstsystem[\sit],\sit)=\superfarth(\ratfrcstsystem[\sit],\sit)\) for all~\(\sit\in\sits\) for which \((\ratfrcstsystem[\sit],\sit)\) is non-degenerate, provided that the map~\(\superfarth(\ratfrcstsystem[\bolleke],\bolleke)\colon\sits\to\reals\) is a positive test supermartingale for~\(\ratfrcstsystem\). \label{lem:prop:two}
\end{enumerate}
\end{lemma}

\begin{proof}
Start from a code for the map~\(\superfarth\colon\preqsits\to\sqgroup{0,+\infty}\) that is lower semicomputable.
By \corollaryref{cor:lsc}, we can invoke a single algorithm that outputs a code \(q\colon\preqsits\times\naturalswithzero\to\rationals\) for~\(\superfarth\) such that \(q(\preqsit,\bolleke)\lsc\superfarth(\preqsit)\) and \(q(\preqsit,n)<q(\preqsit,n+1)\) for all~\(\preqsit\in\preqsits\) and \(n\in\naturalswithzero\).
We'll now use the code \(q\) to construct a code \(q'\) for a lower semicomputable test superfarthingale~\(\superfarth'\in\superfarths\) that satisfies the requirements of the lemma.

Let \(q'\colon\preqsits\times\naturalswithzero\to\rationals\) be defined by~\(q'(\init,n)\coloneqq1\) and
\begin{equation*}
q'(\preqsit\ratfrcst x,n)\coloneqq
\begin{cases}
\max\group[\big]{A(\preqsit,\ratfrcst,x,n)\cup\set{0}}
&\text{if \(\preqsit\ratfrcst x\) is non-degenerate}\\
0 &\text{if \(\preqsit\ratfrcst x\) is degenerate},
\end{cases}
\end{equation*}
for all~\(\preqsit=\preqsitrep\in\preqsits\), \(\ratfrcst\in\ratfrcsts\), \(x\in\posspace\) and \(n\in\naturalswithzero\), where the map
\[
A\colon\smash{\preqsits\times\ratfrcsts\times\posspace\times\naturalswithzero}\to\cset{Q\subseteq\rationals}{\abs{Q}<\infty}
\]
is defined by
\begin{multline}\label{eq:def:A}
A(\preqsit,\ratfrcst,x,n)\coloneqq\cset[\big]{q(\preqsit\ratfrcst x,m)\in\rationals}{0\leq m\leq n\text{, } \\
0\leq q(\preqsit\ratfrcst\andoutcome,m)\text{ and }
\uex_{\ratfrcst}(q(\preqsit\ratfrcst\andoutcome,m))\leq q'(\preqsit,n)}.
\end{multline}
By construction, since the map~\(A\) outputs finite sequences of rationals, the map~\(q'\) is well-defined, non-negative and rational.
It's not too difficult to see that the map~\(A\), and therefore also the map~\(q'\), is recursive.

The map~\(q'\) is non-decreasing in its second argument, as we now show by induction on its first argument.
We start by observing that \(q'(\init,n)\leq q'(\init,n+1)\) for all~\(n\in\naturalswithzero\).
For the induction step, fix any~\(\preqsit=\preqsitrep\in\preqsits\), \(\ratfrcst\in\ratfrcsts\), \(x\in\posspace\) and \(n\in\naturalswithzero\), and assume that \(q'(\preqsit,n)\leq q'(\preqsit,n+1)\).
We then have to show that also \(q'(\preqsit\ratfrcst x,n)\leq q'(\preqsit\ratfrcst x,n+1)\).
This is trivial when \(\preqsit\ratfrcst x\) is degenerate; when \(\preqsit\ratfrcst x\) is non-degenerate, it follows readily from the inequality \(A(\preqsit,\ratfrcst,x,n)\subseteq A(\preqsit,\ratfrcst,x,n+1)\), which is itself immediate from \equationref{eq:def:A}.

For any~\(n\in\naturalswithzero\), the map~\(q'(\bolleke,n)\colon\preqsits\to\reals\) is a test superfarthingale.
To prove this, we may clearly concentrate on the superfarthingale condition.
Fix any~\(\preqsit\in\preqsits\), \(\ratfrcst\in\ratfrcsts\) and \(n\in\naturalswithzero\), and infer from \equationref{eq:def:A} that \(A(\preqsit,\ratfrcst,1,n)=\emptyset\ifandonlyif A(\preqsit,\ratfrcst,0,n)=\emptyset\), so we only need to consider two cases.
If \(A(\preqsit,\ratfrcst,1,n)=A(\preqsit,\ratfrcst,0,n)=\emptyset\), then \(q'(\preqsit\ratfrcst\andoutcome,n)=0\), and therefore trivially \(\uex_{\ratfrcst}(q'(\preqsit\ratfrcst\andoutcome,n))=\uex_{\ratfrcst}(0)=0\leq q'(\preqsit,n)\), where the second equality follows from \ref{axiom:coherence:bounds}.
Otherwise, because the map~\(q\) is increasing in its second argument, there's an \(m\in\set{0,\dots,n}\) such that \(\uex_{\ratfrcst}(q(\preqsit\ratfrcst\andoutcome,m))\leq q'(\preqsit,n)\), with
\[
q(\preqsit\ratfrcst\andoutcome,m)
=\max\group{A(\preqsit,\ratfrcst,\andoutcome,n)\cup\set{0}}
\geq q'(\preqsit\ratfrcst\andoutcome,n),
\]
where the last inequality takes into account that there may be some \(x\in\posspace\) such that \(\preqsit\ratfrcst x\) is degenerate.
Hence, indeed, in this case also
\begin{align*}
\uex_{\ratfrcst}(q'(\preqsit\ratfrcst\andoutcome,n))
\overset{\text{\ref{axiom:coherence:monotonicity}}}{\leq}\uex_{\ratfrcst}(q(\preqsit\ratfrcst\andoutcome,m))
\leq q'(\preqsit,n).
\end{align*}

As a final preliminary step, we infer from \lemmaref{lem:superfarth:bound} that for every (non-degenerate) prequential situation~\(\preqsit\in\preqsits\) there's some real~\(B_\preqsit\in\reals\) such that \(q'(\preqsit,n)\leq B_\preqsit\) for all~\(n\in\naturalswithzero\).

With this set-up phase completed, let \(\superfarth'\) be defined as \(q'(\preqsit,\bolleke)\lsc\superfarth'(\preqsit)\) for all~\(\preqsit\in\preqsits\); note that \(\superfarth'(\init)=1\).
This map is well-defined, real-valued, non-negative and lower semicomputable due to the non-decreasingness, boundedness, non-negativity and recursiveness of \(q'\) respectively, so we only need to check the superfarthingale property explicitly in order to conclude that \(\superfarth'\) is a lower semicomputable test superfarthingale.
To this end, fix any~\(\preqsit\in\preqsits\) and \(\ratfrcst\in\ratfrcsts\).
If we recall that the map~\(q(\bolleke,n)\colon\preqsits\to\reals\) is a test superfarthingale for every \(n\in\naturalswithzero\), we immediately infer from \ref{axiom:coherence:uniform:convergence} and the real-valuedness of \(\superfarth'\) that \(\uex_{\ratfrcst}(\superfarth'(\preqsit\ratfrcst\andoutcome))=\lim_{n\to\infty}\uex_{\ratfrcst}(q'(\preqsit\ratfrcst\andoutcome,n))\leq\lim_{n\to\infty}q'(\preqsit,n)=\superfarth'(\preqsit)\).

We are done if we can show that \(F'\) satisfies the conditions~\ref{lem:prop:one} and~\ref{lem:prop:two}.
For \ref{lem:prop:one}, fix any degenerate prequential situation~\(\preqsit\in\preqsits\) and note that then \(q'(\preqsit,n)=0\) for all~\(n\in\naturalswithzero\) by construction. Hence, indeed, \(\superfarth'(\preqsit)=0\).

For \ref{lem:prop:two}, fix any rational forecasting system~\(\ratfrcstsystem\in\ratfrcstsystems\), consider the map~\(\test\colon\sits\to\reals\) defined by~\(\test(\sit)\coloneqq\superfarth(\ratfrcstsystem[\sit],\sit)\) for all~\(\sit\in\sits\), and assume that \(\test\) is a positive test supermartingale.
We must now show that \(\superfarth'(\ratfrcstsystem[\sit],\sit)=\test(\sit)\) for all~\(\sit\in\sits\) for which the prequential situation~\((\ratfrcstsystem[\sit],\sit)\) is non-degenerate.

By construction, \(\superfarth'(\ratfrcstsystem[\sit],\sit)\leq\superfarth(\ratfrcstsystem[\sit],\sit)=\test(\sit)\) for all~\(\sit\in\sits\).
Assume towards contradiction that there's some \(\altsit\in\sits\) for which \((\ratfrcstsystem[\altsit],\altsit)\) is non-degenerate and \(\superfarth'(\ratfrcstsystem[\altsit],\altsit)<\test(\altsit)\), implying that there's some \(\epsilon>0\) such that \(q'((\ratfrcstsystem[\altsit],\altsit),n)+\epsilon<\test(\altsit)\) for all~\(n\in\naturalswithzero\).
We'll use an induction argument to show that this is impossible.

Since, by assumption,
\[
q((\ratfrcstsystem[\altsit],\altsit),\bolleke)\lsc\test(\altsit)>0
\text{ and }
q((\ratfrcstsystem[\altsit],\altsit),n)<q((\ratfrcstsystem[\altsit],\altsit),n+1)
\text{ for all~\(n\in\naturalswithzero\)},
\]
there are \(\epsilon_0,\epsilon_1,\dots,\epsilon_{\abs{\altsit}}\in\reals\) and \(n_0,n_1,\dots,n_{\abs{\altsit}}\in\naturalswithzero\) such that
\begin{align}
0<\epsilon_0<\epsilon_1<\dots<\epsilon_{\abs{\altsit}}<\epsilon\label{eq:eps}\\
\test(\altsitto{\ell})<q((\ratfrcstsystem[\altsitto{\ell}],\altsitto{\ell}),n_\ell)+\epsilon_\ell\label{eq:app:above}\\
0\leq q((\ratfrcstsystem[\altsitto{k}]\ratfrcstsystem(\altsitto{k}),\altsitto{k}\andoutcome),n_{k+1})\label{eq:positivity}\\
q((\ratfrcstsystem[\altsitto{k}]\ratfrcstsystem(\altsitto{k}),\altsitto{k}\andoutcome),n_{k+1})+\epsilon_{k}<\test(\altsitto{k}\andoutcome)\label{eq:app:below}
\end{align}
for all~\(k\in\set{0,1,\dots,\abs{\altsit}-1}\) and \(\ell\in\set{0,1,\dots,\abs{\altsit}}\).
The argument starts with~\(\ell\coloneqq\abs{\altsit}\) and \(k\coloneqq\abs{\altsit}-1\), finding \(\epsilon_\ell\) such that \eqref{eq:eps} is satisfied, and finding \(n_{k+1}\) such that \eqref{eq:app:above} and \eqref{eq:positivity} are satisfied.
We then move to \(\ell\coloneqq\abs{\altsit}-1\) and \(k\coloneqq\abs{\altsit}-2\), find \(\epsilon_\ell\) such that \eqref{eq:eps} and \eqref{eq:app:below} are satisfied, and find \(n_{k+1}\) such that \eqref{eq:app:above} and \eqref{eq:positivity} are satisfied.
And so on \dots; these conditions are depicted below for a situation~\(\altsit\in\sits\) for which \(\abs{\altsit}=5\).

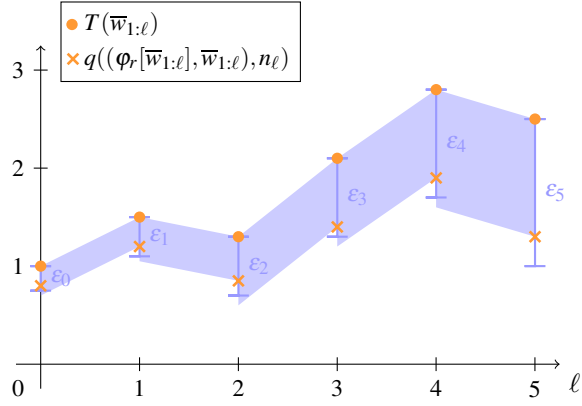
\begin{center}
\begin{tikzpicture}[scale=1.3,
  orangenode/.style={circle,fill,orange!80!white,inner sep=1.5pt},
  orangecross/.style={cross out,draw,orange!80!white,inner sep=1.5pt,line width=1pt}
]

\draw[->] (-0.25,0) -- (5.25,0) node[below right] {\(\ell\)};
\draw[->] (0,-0.25) -- (0,3.25); 

\foreach \x in {1,2,3,4,5}
    \draw (\x,2pt) -- (\x,-2pt)
	node[anchor=north] {\small\(\x\)};
\foreach \y in {1,2,3}
    \draw (2pt,\y) -- (-2pt,\y) 
    node[anchor=east] {\small\(\y\)};
\node[anchor=north east] at (-2pt, -2pt)  {\small\(0\)};

\def\ya{1} 
\def\yb{1.5}
\def\yc{1.3}
\def\yd{2.1}
\def\ye{2.8}
\def\yf{2.5}
\coordinate (A) at (0,\ya); 
\coordinate (B) at (1,\yb);
\coordinate (C) at (2,\yc);
\coordinate (D) at (3,\yd);
\coordinate (E) at (4,\ye);
\coordinate (F) at (5,\yf);
\def\a{0.25} 
\def\b{0.4}
\def\c{0.6}
\def\d{0.8}
\def\e{1.1}
\def\f{1.5}
\def\width{3pt}
\def\qa{0.05} 
\def\qb{0.1}
\def\qc{0.15}
\def\qd{0.1}
\def\qe{0.2}
\def\qf{0.3}
\coordinate (QA) at (0,\ya-\a+\qa); 
\coordinate (QB) at (1,\yb-\b+\qb);
\coordinate (QC) at (2,\yc-\c+\qc);
\coordinate (QD) at (3,\yd-\d+\qd);
\coordinate (QE) at (4,\ye-\e+\qe);
\coordinate (QF) at (5,\yf-\f+\qf);

\draw[line width=0.8pt,blue!40!white] (A) --+ (-\width,0) --+ (\width,0) --+ (0pt,0) --++ (0,-\a) node [midway, anchor=west] {\small\(\epsilon_0\)} --+ (-\width,0) --+ (\width,0);
\draw[line width=0.8pt,blue!40!white] (B) --+ (-\width,0) --+ (\width,0) --+ (0pt,0) --++ (0,-\b) node [midway, anchor=west] {\small\(\epsilon_1\)} --+ (-\width,0) --+ (\width,0);
\draw[line width=0.8pt,blue!40!white] (C) --+ (-\width,0) --+ (\width,0) --+ (0pt,0) --++ (0,-\c) node [midway, anchor=west] {\small\(\epsilon_2\)} --+ (-\width,0) --+ (\width,0);
\draw[line width=0.8pt,blue!40!white] (D) --+ (-\width,0) --+ (\width,0) --+ (0pt,0) --++ (0,-\d) node [midway, anchor=west] {\small\(\epsilon_3\)} --+ (-\width,0) --+ (\width,0);
\draw[line width=0.8pt,blue!40!white] (E) --+ (-\width,0) --+ (\width,0) --+ (0pt,0) --++ (0,-\e) node [midway, anchor=west] {\small\(\epsilon_4\)} --+ (-\width,0) --+ (\width,0);
\draw[line width=0.8pt,blue!40!white] (F) --+ (-\width,0) --+ (\width,0) --+ (0pt,0) --++ (0,-\f) node [midway, anchor=west] {\small\(\epsilon_5\)} --+ (-\width,0) --+ (\width,0);

\fill[opacity=0.5,blue!40!white] (A) -- (B) -- (QB) -- (0,\ya-\b+\qb) --cycle;
\fill[opacity=0.5,blue!40!white] (B) -- (C) -- (QC) -- (1,\yb-\c+\qc) --cycle;
\fill[opacity=0.5,blue!40!white] (C) -- (D) -- (QD) -- (2,\yc-\d+\qd) --cycle;
\fill[opacity=0.5,blue!40!white] (D) -- (E) -- (QE) -- (3,\yd-\e+\qe) --cycle;
\fill[opacity=0.5,blue!40!white] (E) -- (F) -- (QF) -- (4,\ye-\f+\qf) --cycle;

\foreach \position in {A,B,C,D,E,F}
    \node [orangenode] at (\position) {};
\foreach \position in {QA,QB,QC,QD,QE,QF}
    \node [orangecross] at (\position) {};

\matrix [draw,below right,inner sep=2pt] at (0.2,3.7) {
  \node [orangenode,label=right:\small\(\test(\altsitto{\ell})\)] {}; \\
  \node [orangecross,label=right:\small{\(q((\ratfrcstsystem[\altsitto{\ell}],\altsitto{\ell}),n_\ell)\)}] {}; \\
};

\end{tikzpicture}

\end{center}
\noindent Now, let \(N\coloneqq\max\set{n_0,n_1,\dots,n_{\abs{\altsit}}}\).
To start the induction argument, observe that, trivially, \(q'(\init,N)=1>\test(\init)-\epsilon_0\).
For the induction step, we fix any~\(k\in\set{0,1,\dots,\abs{\altsit}-1}\) and assume that \(q'((\ratfrcstsystem[\altsitto{k}],\altsitto{k}),N)>\test(\altsitto{k})-\epsilon_k\).
It then follows that
\begin{align*}
\uex_{\ratfrcstsystem(\altsitto{k})}\group[\big]{q((\ratfrcstsystem[\altsitto{k}]\ratfrcstsystem(\altsitto{k}),\altsitto{k}\andoutcome),n_{k+1})}
\overset{\eqref{eq:app:below},\text{\ref{axiom:coherence:monotonicity}}}&{\leq}
\uex_{\ratfrcstsystem(\altsitto{k})}(\test(\altsitto{k}\andoutcome)-\epsilon_{k}) \\
\overset{\text{\ref{axiom:coherence:constantadditivity}}}&{=}
\uex_{\ratfrcstsystem(\altsitto{k})}(\test(\altsitto{k}\andoutcome))-\epsilon_{k} \\
&\leq\test(\altsitto{k})-\epsilon_{k} \\
&\leq q'((\ratfrcstsystem[\altsitto{k}],\altsitto{k}),N),
\end{align*}
where the penultimate inequality follows form the assumption that \(\test\) is a supermartingale, and the last inequality from the induction hypothesis.
Hence, by Equations~\eqref{eq:def:A} and~\eqref{eq:positivity},
\begin{equation*}
q((\ratfrcstsystem[\altsitto{k+1}],\altsitto{k+1}),n_{k+1})
\in
A((\ratfrcstsystem[\altsitto{k}],\altsitto{k}),\ratfrcstsystem(\altsitto{k}),\altsitat{k+1},N),
\end{equation*}
which implies that
\begin{align*}
q'((\ratfrcstsystem[\altsitto{k+1}],\altsitto{k+1}),N)
&\geq\max A((\ratfrcstsystem[\altsitto{k}],\altsitto{k}),\ratfrcstsystem(\altsitto{k}),\altsitat{k+1},N) \\
&\geq q((\ratfrcstsystem[\altsitto{k+1}],\altsitto{k+1}),n_{k+1}) \\
\overset{\eqref{eq:app:above}}&{>}
\test(\altsitto{k+1})-\epsilon_{k+1}.
\end{align*}
Repeating this argument until we reach \(k=\abs{\altsit}-1\), we eventually find that \(q'((\ratfrcstsystem[\altsit],\altsit),N)\allowbreak>\test(\altsit)-\epsilon_{\abs{\altsit}}>\test(\altsit)-\epsilon\), which is the desired contradiction.
\end{proof}
\noindent The following result is now immediate.

\begin{corollary} \label{cor:uniform:program}
There's a single algorithm that, upon the input of a code for a lower semicomputable map~\(\superfarth\colon\preqsits\to\sqgroup{0,+\infty}\), outputs a code for a lower semicomputable test superfarthingale~\(\superfarth'\in\superfarths\) such that, for all prequential situations~\(\preqsit\in\preqsits\),
\begin{enumerate}[label=\upshape(\roman*),leftmargin=*,itemsep=0pt]
\item \(\superfarth'(\preqsit)=0\) if \(\preqsit\) is degenerate;
\item \(\superfarth'(\preqsit)=\superfarth(\preqsit)\) if \(\preqsit\) is non-degenerate and \(\superfarth\) is a positive test superfarthingale. 
\end{enumerate}
\end{corollary}
\end{document}